\documentclass[12pt]{amsart}

 \usepackage{amsfonts,graphics,amsmath,amsthm,amsfonts,amscd,amssymb,amsmath,latexsym,multicol,
 mathrsfs}
\usepackage{epsfig,url}
\usepackage{flafter}
\usepackage{fancyhdr}
\usepackage{hyperref}
\hypersetup{colorlinks=true, linkcolor=black}

\makeatletter

\def\jobis#1{FF\fi
  \def\predicate{#1}%
  \edef\predicate{\expandafter\strip@prefix\meaning\predicate}%
  \edef\job{\jobname}%
  \ifx\job\predicate
}

\makeatother

\if\jobis{proposal}%
\else
\fi

 \usepackage[matrix, arrow]{xy}

\DeclareMathOperator{\Supp}{Supp}

\DeclareMathOperator{\Pic}{Pic}

\DeclareMathOperator{\Div}{Div}

 \numberwithin{equation}{subsection}
 \numberwithin{footnote}{subsection}

 \newtheorem{lem}[subsection]{Lemma}
 \newtheorem{prop}[subsection]{Proposition}
 \newtheorem{thm}[subsection]{Theorem}

{
    \newtheoremstyle{upright}%
        {8pt plus2pt minus4pt}%
        {8pt plus2pt minus4pt}%
        {\upshape}%
        {}%
        {\bfseries\scshape}%
        {}%
        {1em}%
        {}%
\theoremstyle{upright}

 \newtheorem{defn}[subsection]{Definition}

 \newtheorem{rem}[subsection]{Remark}

}

 \newcommand{\N}{\mathbb N}
 \newcommand{\PP}{\mathbb P}
 
 \newcommand{\Q}{\mathbb Q}
 \newcommand{\R}{\mathbb R}
 \newcommand{\E}{\mathbb E}
 \newcommand{\Z}{\mathbb Z}

 \newcommand{\bir}{\dashrightarrow}
 \newcommand{\lin}{\sim}
 \newcommand{\num}{\equiv}
 
 \newcommand{\rddown}[1]{\left\lfloor{#1}\right\rfloor} 

\title{\large E\MakeLowercase{xistence of} M\MakeLowercase{ori fibre spaces for 3-folds in char $p$}}
\thanks{2010 MSC: 14E30}
\author{Caucher Birkar and Joe Waldron}
\date{\today}
\begin{document}
\maketitle

\begin{abstract} 
We prove the following results for projective klt pairs of dimension $3$ over an algebraically 
closed field of char $p>5$: the cone theorem, the base point free theorem, 
the contraction theorem, finiteness of minimal models, termination with scaling, existence of 
Mori fibre spaces, etc. 
\end{abstract}

\tableofcontents


\section{Introduction}

We work over an algebraically closed field $k$ (mostly of char $p>5$). Boundary divisors are always 
assumed to be with real coefficients unless otherwise stated. 

After [\ref{HX}], many results concerning the log minimal model program (LMMP) for $3$-folds over $k$ of char $p>5$ 
were settled in [\ref{B}] such as existence of log flips and log minimal models, special cases of the 
base point free and contraction theorems, special cases of Koll\'ar-Shokurov connectedness 
principle, existence of $\Q$-factorial dlt models, ACC for log canonical thresholds, etc. 
One of the main problems not treated in [\ref{B}] is the existence of Mori fibre spaces.  
Their existence is proved in this paper. We also settle various other problems that are 
discussed below.\\

{\textbf{\sffamily{Cone theorem.}}}
In characteristic $0$, the base point free theorem and the cone and contraction 
theorems are among the first results of the LMMP. They are derived from  
the Kawamata-Viehweg vanishing theorem. 
But the story in positive characteristic is quite different because such vanishing theorems fail. 
Unlike in characteristic $0$, existence of flips and minimal models is a fundamental ingredient of the 
proof of the cone and contraction and base point free theorems.

\begin{thm}\label{t-cone}
Let $(X,B)$ be a projective $\Q$-factorial dlt pair of dimension $3$ over $k$ of char $p>5$. 
Then there is a countable number of rational curves $\Gamma_i$ such that\\ 

$(i)$ $\overline{NE}(X)=\overline{NE}(X)_{K_X+B\ge 0}+\sum_i \R [\Gamma_i]$, 

$(ii)$ $-6\le (K_X+B)\cdot \Gamma_i<0$, 

$(iii)$ for any ample $\R$-divisor $A$, 
$$
(K_X+B+A)\cdot \Gamma_i\ge 0
$$ 
for all but finitely many $i$, and 

$(iv)$ the rays $\R [\Gamma_i]$ do not accumulate inside  $\overline{NE}(X)_{K_X+B<0}$.\\
\end{thm}

The theorem is proved in Section \ref{s-cone} where we also prove some other results concerning 
extremal rays (see Proposition \ref{p-polytope-rays}).
Special cases of the theorem were proved in [\ref{Keel}, Proposition 0.6 ][\ref{CTX}, Theorem 1.7].\\

{\textbf{\sffamily{Base point freeness.}}}
The proof of the next result, given in Section \ref{s-main-results}, relies 
on the results in Sections \ref{s-cone} to \ref{s-bpf}.
The whole proof occupies a big chunk of this paper and it contains many of our key ideas and 
technical results.

\begin{thm}\label{t-bpf}
Let $(X,B)$ be a projective klt pair of dimension $3$ over $k$ of char $p>5$ and $X\to Z$ a 
projective contraction. 
Assume that $D$ is an $\R$-divisor such that $D$ is nef$/Z$ and $D-(K_X+B)$ is nef and big$/Z$. Then 
 $D$ is  semi-ample$/Z$.
\end{thm}

The theorem was proved in [\ref{B}][\ref{Xu}] when $D$ is a big $\Q$-divisor using existence of 
minimal models and Keel's semi-ampleness techniques. When $D$ is not big, Keel's methods do not 
apply, at least not directly. To deal with this issue, in [\ref{CTX}], a canonical bundle formula is used to 
reduce the problem to surfaces when $D$ is a $\Q$-divisor with numerical dimension $\nu(D)=2$ 
and assuming that $B+A$ is a $\Q$-boundary with coefficients $>\frac{2}{p}$ for some $0\le A\sim_\R D-(K_X+B)$. 
The canonical bundle formula is derived from the 
theory of moduli of pointed curves. 
Our proof is very different and it does not involve canonical bundle formulas.\\

{\textbf{\sffamily{Contraction theorem.}}}
The next result is a consequence of the base point free theorem and the cone theorem above.

\begin{thm}\label{t-contraction}
Let $(X,B)$ be a projective $\Q$-factorial dlt pair of dimension $3$ over $k$ of char $p>5$ 
and $X\to Z$ a projective contraction. 
Then any $K_X+B$-negative extremal ray$/Z$ can be contracted by a projective contraction.
\end{thm}

The proof is given in Section \ref{s-main-results}.
The theorem was proved in [\ref{B}, Theorem 1.5][\ref{Xu}] for extremal rays of 
flipping or divisorial type.\\

{\textbf{\sffamily{Finiteness of minimal models.}}}
We prove finiteness of minimal models under suitable assumptions and 
derive termination with scaling. This is similar to the characteristic $0$ case [\ref{BCHM}].

First we introduce some notation. Let $X\to Z$ be a projective contraction of normal projective varieties  
over $k$ of char $p>5$ where $X$ is $\Q$-factorial of dimension $3$. Let $A\ge 0$ be a $\Q$-divisor on $X$, and 
$V$ a rational finite dimensional affine subspace of the space of $\R$-Weil divisors on $X$. Define 
$$
\mathcal{L}_{A}(V)=\{\Delta \mid 0\le (\Delta-A)\in V, ~ \mbox{and $(X,\Delta)$ is lc} \}
$$
As in [\ref{Shokurov}, 1.3.2], one can show that $\mathcal{L}_{A}(V)$ 
is a rational polytope inside the rational affine space $A+V$, that is, 
it is the convex hull of finitely many rational points in $A+V$: 
this follows from existence of log resolutions.

\begin{thm}\label{t-finiteness}
Under the above setting, assume in addition that $A$ is big$/Z$. Let $\mathcal{C}\subseteq \mathcal{L}_{A}(V)$ 
be a rational polytope such that $(X,\Delta)$ is klt for every $\Delta\in \mathcal{C}$. Then there are finitely many birational 
maps $\phi_i\colon X\bir Y_i/Z$ such that for any $\Delta\in \mathcal{C}$ with $K_X+\Delta$ pseudo-effective$/Z$, there 
is $i$ such that $(Y_i,\Delta_{Y_i})$ is a log minimal model of $(X,\Delta)$ over $Z$.
\end{thm}

As usual $\Delta_{Y_i}$ means the pushdown $(\phi_i)_*\Delta$. A conditional proof of the theorem 
is given in Section \ref{s-finiteness}. At the end in Section \ref{s-main-results} the extra assumptions are removed. \\ 

{\textbf{\sffamily{Termination with scaling.}}}
Minimal models were constructed in [\ref{B}] by a rather indirect approach. 
It is useful in many situations to know that running an LMMP ends up with a minimal model. 
It is even more important for constructing Mori fibre spaces. 

\begin{thm}\label{t-term}
Let $(X,B+C)$ be a $\Q$-factorial projective klt pair of dimension $3$ over $k$ of char $p>5$ 
and $X\to Z$ a projective contraction. 
Assume that $B\ge 0$ is big$/Z$,  
$C\ge 0$ is $\R$-Cartier, and $K_X+B+C$ is nef$/Z$. Then we can run the LMMP$/Z$ on $K_X+B$ with scaling of $C$ 
and it terminates.
\end{thm}

\begin{thm}\label{t-term-2}
Let $(X,B+C)$ be a $\Q$-factorial projective klt pair of dimension $3$ over $k$ of char $p>5$ 
and $X\to Z$ a projective contraction. 
Assume $C\ge 0$ is ample, and $K_X+B+C$ is nef$/Z$. Then we can run the LMMP$/Z$ on $K_X+B$ with scaling of $C$ 
and it terminates.
\end{thm}

The proofs are given in Section \ref{s-finiteness} under certain assumptions. Unconditional proofs 
are in Section \ref{s-main-results}.\\

{\textbf{\sffamily{Mori fibre spaces.}}}
Finally we come to the result which is the title of this paper. 

\begin{thm}\label{t-Mfs}
Let $(X,B)$ be a dlt pair of dimension $3$ over $k$ of char $p>5$ and $X\to Z$ a projective 
contraction. Assume $K_X+B$ is not pseudo-effective$/Z$. Then $(X,B)$ has a Mori fibre space$/Z$.
If $X$ is $\Q$-factorial, then we can run an LMMP$/Z$ on $K_X+B$ which ends with a Mori fibre space$/Z$. 
\end{thm}

The proof is given at the very end of the paper in Section \ref{s-main-results}.
This theorem combined with [\ref{B}, Theorem 1.2] say that any  klt pair $(X,B)$ of 
dimension $3$ over $k$ of char $p>5$, projective over some base $Z$, either has a log minimal model 
or a Mori fibre space over $Z$.\\

{\textbf{\sffamily{Acknowledgements.}}}
This work was partially supported by a grant of the Leverhulme Trust.
Part of this work was done when the first author visited National Taiwan University in August-September 2014 with the support of the 
Mathematics Division (Taipei Office) of the National Center for Theoretical Sciences, and the visit was arranged by Jungkai A. Chen. 
He wishes to thank them for their hospitality.\\

\section{Preliminaries}

All the varieties and algebraic spaces in this paper are defined over an algebraically closed field 
$k$ unless otherwise stated. 

\subsection{Contractions and divisors endowed with a map}\label{ss-contractions} 
A \emph{contraction} $f\colon X\to Z$ of algebraic spaces over $k$ is a proper morphism such that 
$f_*\mathcal{O}_X=\mathcal{O}_Z$. 
When  $X,Z$ are quasi-projective varieties over $k$ and $f$ is 
projective, we refer to $f$ as a \emph{projective contraction} to avoid confusion.
In this case, by a \emph{fibre} of $f$ we always mean a scheme-theoretic 
fibre unless stated otherwise.

For a nef $\Q$-divisor $L$ on a projective scheme $X$ over $k$, the 
\emph{exceptional locus} $\mathbb{E}(L)$ is the union of those 
positive-dimensional integral subschemes $Y\subseteq X$ such that $L|_Y$ is not big, i.e. $(L|_Y)^{\dim Y}=0$. 
We say $L$ is \emph{endowed with a map} $f\colon X\to V$, where $V$ is an algebraic space over $k$ and  
$f$ is a proper morphism, if an integral subscheme $Y$ is contracted by $f$ (i.e. $\dim Y>\dim f(Y)$) 
if and only if $L|_Y$ is not big.

\subsection{Rational maps}

\begin{lem}\label{l-nef-pushdown}
Let $\phi\colon X\bir Y$ be a birational map between normal projective varieties over $k$. 
Assume $\phi^{-1}$ does not contract divisors. Let $D$ be a nef $\R$-divisor on $X$ 
such that $D_Y=\phi_*D$ is $\R$-Cartier. Let $f\colon W\to X$ and $g\colon W\to Y$ 
be a common resolution. Then, $E:=g^*D_Y-f^*D$ is effective and exceptional$/Y$.
\end{lem}
\begin{proof}
It is obvious that $E$ is exceptional$/Y$. The effectivity is a consequence of the 
negativity lemma.\\
\end{proof}

\begin{lem}\label{l-open-isom}
Let $\phi\colon X\bir Y$ be a birational map between normal projective varieties over $k$. 
Assume $\phi^{-1}$ does not contract divisors. Then there is an open subset $U\subseteq X$ 
such that $\phi|_U$ is an isomorphism and codimension of $Y\setminus \phi(U)$ 
is at least $2$. 
\end{lem} 
\begin{proof}
Let $U\subseteq X$ be the largest open set such that $\phi|_U$ is an isomorphism.
Assume that codimension of $Y\setminus \phi(U)$ is one and let $S$ be one of its components 
of codimension one. Since $\phi^{-1}$ does not contract divisors, $\phi^{-1}$ is an isomorphism 
near the generic point of $S$. This means that there is an open set $V\subseteq X$ intersecting the 
birational transform of $S$ such that $\phi|_V$ is an isomorphism. But then 
$\phi|_{U\cup V}$ is an isomorphism which contradicts the maximality of $U$.\\  
\end{proof}

\subsection{Pairs}\label{ss-pairs} 
A \emph{pair} $(X,B)$ consists of a normal quasi-projective variety $X$ over $k$  
and an \emph{$\R$-boundary} $B$, that is an $\R$-divisor $B$ on $X$ with coefficients in $[0,1]$, 
such that $K_X+B$ is $\mathbb{R}$-Cartier. When $B$ has rational coefficients we 
say $B$ is a \emph{$\Q$-boundary}.
We say that $(X,B)$ is \emph{log smooth} if $X$ is smooth 
and $\Supp B$ has simple normal crossing singularities.

Let $(X,B)$ be a pair. For a prime divisor $D$ on some birational model of $X$ with a
nonempty centre on $X$, $a(D,X,B)$ denotes the \emph{log discrepancy} which is defined 
by taking a projective birational morphism $f\colon Y\to X$ from a normal variety 
containing $D$ as a prime divisor and putting $a(D,X,B)=1-b$ where $b$ is the 
coefficient of $D$ in $B_Y$ and $K_Y+B_Y=f^*(K_X+B)$.  

 As in characteristic $0$, we can define various types of singularities using 
log discrepancies. Let $(X,B)$ be a pair.
We say that the pair is \emph{log canonical} or lc for short (resp. \emph{Kawamata log terminal} or {klt} 
for short) if 
$a(D,X,B)\ge 0$ (resp. $a(D,X,B)>0$) for any prime divisor $D$ on birational models of $X$. 
An \emph{lc centre} of $(X,B)$ is the image in $X$ of a $D$ with $a(D,X,B)=0$.
On the other hand, we say that $(X,B)$ is \emph{dlt} if there is a closed subset $P\subset X$ 
such that $(X,B)$ is log smooth 
outside $P$ and no lc centre of $(X,B)$ is inside $P$. In particular, 
the lc centres of $(X,B)$ are exactly the components of $S_1\cap \cdots \cap S_r$ where $S_i$ 
are among the components of $\rddown{B}$. Moreover, there is a log resolution 
$f\colon Y\to X$ of $(X,B)$ such that $a(D,X,B)>0$ for any prime divisor $D$ on $Y$ 
which is exceptional$/X$, eg take a log resolution $f$ 
which is an isomorphism over $X\setminus P$.  
Finally, we say that $(X,B)$ is \emph{plt} if it is dlt and each connected component of 
$\rddown{B}$ is irreducible. In particular, the only lc centres of $(X,B)$ are 
the components of $\rddown{B}$.

\subsection{Minimal models and Mori fibre spaces}\label{ss-mmodels} 
Let $(X,B)$ be a pair and $X\to Z$ a projective contraction over $k$.
A pair $(Y,B_Y)$ with a projective contraction $Y\to Z$ and a birational map
$\phi\colon X\bir Y/Z$ is a \emph{log birational model} of $(X,B)$ 
if  $B_Y$ is the sum of the birational transform of $B$ 
and the reduced exceptional divisor of $\phi^{-1}$. 
We say that $(Y,B_Y)$ is a \emph{weak lc model} of $(X,B)$ over $Z$ if in addition\\\\
(1) $K_Y+B_Y$ is nef/$Z$, and\\
(2) for any prime divisor $D$ on $X$ which is exceptional/$Y$, we have
$$
a(D,X,B)\le a(D,Y,B_Y)
$$

And  we call $(Y,B_Y)$ a \emph{log minimal model} of $(X,B)$ over $Z$ if in addition\\\\
(3) $(Y,B_Y)$ is $\Q$-factorial dlt, and\\
(4) the inequality in (2) is strict.\\

A weak lc model or log minimal model $(Y,B_Y)$ is said to be \emph{good} if $K_Y+B_Y$ is semi-ample$/Z$.
When $K_X+B$ is big$/Z$, the \emph{lc model} of $(X,B)$ over $Z$ is a weak lc model $(Y,B_Y)$
over $Z$ with $K_Y+B_Y$ ample$/Z$.

On the other hand, a log birational model $(Y,B_Y)$ of $(X,B)$ is called a 
\emph{Mori fibre space} of $(X,B)$ over $Z$ if there is a $K_Y+B_Y$-negative extremal projective 
contraction $Y\to T/Z$, and if for any prime divisor $D$ on birational models of $X$ we have
$$
a(D,X,B)\le a(D,Y,B_Y)
$$
with strict inequality if $D\subset X$ and if it is exceptional/$Y$,

Note that the above definitions are the same as in [\ref{B}] but 
slightly different from the traditional definitions in that we allow $\phi^{-1}$ to contract divisors. However, 
if $(X,B)$ is plt (hence also klt) the definitions coincide. Actually in this paper we usually 
deal with models such that $\phi^{-1}$ does not contract divisors.

Let $(X,B)$ be an lc pair over $k$. A $\Q$-factorial dlt pair $(Y,B_Y)$ is 
a \emph{$\Q$-factorial dlt model} of $(X,B)$ if there is a projective birational 
morphism $f\colon Y\to X$ such that $K_Y+B_Y=f^*(K_X+B)$ and such that every exceptional prime
divisor of $f$ has coefficient $1$ in $B_Y$.  

One of the fundamental ingredients of this paper is the following result which was 
proved after the developments in [\ref{HX}].

\begin{thm}[{[\ref{B}]}]\label{t-mmodel}
Let $(X,B)$ be a klt pair of dimension $3$ over $k$ of char $p>5$ 
and $X\to Z$ a projection contraction. If  
$K_X+B$ is pseudo-effective over $Z$, then $(X,B)$ has a log minimal model over $Z$.  
\end{thm}

\subsection{LMMP with scaling}\label{s-scaling}
Let $(X,B+C)$ be an lc pair and $X\to Z$ a projective contraction over $k$ such that 
$K_X+B+C$ is nef/$Z$, $B\ge 0$, and $C\ge 0$ is $\R$-Cartier. 
Suppose that either $K_X+B$ is nef/$Z$ or there is an extremal ray $R/Z$ such
that $(K_X+B)\cdot R<0$ and $(K_X+B+\lambda_1 C)\cdot R=0$ where
$$
\lambda_1:=\inf \{t\ge 0~|~K_X+B+tC~~\mbox{is nef/$Z$}\}
$$
 If $R$ defines a Mori fibre structure, we stop. Otherwise assume that $R$ gives a divisorial 
contraction or a log flip $X\bir X'$. We can now consider $(X',B'+\lambda_1 C')$  where $B'+\lambda_1 C'$ is 
the birational transform 
of $B+\lambda_1 C$ and proceed similarly. That is, suppose that either $K_{X'}+B'$ is nef/$Z$ or 
there is an extremal ray $R'/Z$ such
that $(K_{X'}+B')\cdot R'<0$ and $(K_{X'}+B'+\lambda_2 C')\cdot R'=0$ where
$$
\lambda_2:=\inf \{t\ge 0~|~K_{X'}+B'+tC'~~\mbox{is nef/$Z$}\}
$$
and so on. By continuing this process, we obtain a 
special kind of LMMP$/Z$ which is called the \emph{LMMP$/Z$ on $K_X+B$ with scaling of $C$}; note that it is not unique. 
When we speak of running such an LMMP we make sure that all the necessary ingredients exist, eg the 
extremal rays, the contractions of the rays, etc.\\

\subsection{Nef reduction maps}

\begin{thm}[{[\ref{BCE}][\ref{CTX}]}]\label{t-nef-reduction}
Let $X$ be a normal projective variety over an uncountable $k$, and $L$ a nef $\R$-divisor on $X$. 
Then there is a rational map $f\colon X\bir Z$ and a nonempty open subset $V\subseteq Z$ such that 

$\bullet$ $f$ is proper over $V$,

$\bullet$ $L|_F\equiv 0$ for the very general fibres $F$ of $f$ over $V$, and 

$\bullet$ if $x\in X$ is a very general point and $C$ a curve passing through $x$, then 
$L\cdot C=0$ iff $C$ is inside the fibre of $f$ containing $x$.    
\end{thm}

We call $f$ a \emph{nef reduction map} of $L$ and call $\dim Z$ the \emph{nef dimension} of $L$ 
and denote it by $n(L)$. If $k$ is countable we can define $n(L)$ to be the nef dimension of $L$ 
after extending $k$ to an uncountable 
algebraically closed field.

The theorem was proved in [\ref{BCE}] for $k$ of characteristic $0$ and $L$ a $\Q$-divisor. It was remarked in 
[\ref{CTX}, 2.4] that the proof in [\ref{BCE}] also works in char $p>0$ and for $L$ an $\R$-divisor. 

In general nef divisors with maximal nef dimension are far from being big or even with 
nonnegative Kodaira dimension. However, nef log divisors behave much better in this sense 
as the following statement shows.

\begin{thm}[{[\ref{CTX}][\ref{B}]}]\label{t-nef-dimension}
Let $(X,B)$ be a projective pair over $k$ such that $B$ is big and $K_X+B$ is nef. 
If the nef dimension $n(K_X+B)=\dim X$, then $K_X+B$ is big.
\end{thm}

The theorem was proved in [\ref{CTX}, Theorem 1.4] in char $p>0$. A short proof was given in [\ref{B}, Theorem 1.11] 
in any characteristic. The theorem actually also holds if $B$ is not a boundary, i.e. has arbitrary coefficients.

\section{Extremal rays and the cone theorem}\label{s-cone}

\subsection{The cone theorem} 
In this subsection we will prove Theorem \ref{t-cone}. First we do some preparations.

\begin{lem}\label{l-f-d-ray}
Let $(X,B)$ be a $\Q$-factorial projective dlt pair of dimension $3$ over 
$k$ of char $p>5$. Suppose that $R$ is a $K_X+B$-negative extremal ray such that 
$N\cdot R=0$ for some nef and big $\Q$-Cartier divisor $N$.
 Then $R$ is generated by some  rational curve 
$\Gamma$ with $(K_X+B)\cdot \Gamma\ge -3$. 
\end{lem}
\begin{proof}
Perturbing the coefficients of $B$ we can assume $(X,B)$ is klt and $B$ is a $\Q$-divisor. 
Since $N\cdot R=0$ for some nef and big $\Q$-Cartier divisor $N$, by [\ref{B}, 3.3], 
there is an ample $\Q$-divisor $A$ such that 
$L=K_X+B+A$ is nef and big and $L^\perp=R$. Moreover, by [\ref{B}, 1.4 and 1.5], 
$L$ is semi-ample and $R$ can be contracted via a projective birational contraction $X\to Z$ 
which is either a flipping or a divisorial contraction.

First suppose $X\to Z$ is a divisorial contraction and let $S$ be the  contracted divisor. 
Let $b$ be the coefficient of $S$ in $B$ and let $\Delta=B+(1-b)S$. 
By adjunction we can write 
$$
K_{S^\nu}+\Delta_{S^\nu}=(K_X+\Delta)|_{S^\nu}
$$
where ${S^\nu}$ is the normalization of $S$ and $\Delta_{S^\nu}\ge 0$ [\ref{Keel}, 5.3]. 
Let $S'\to {S^\nu}$ be the minimal resolution and $S'\to V$ the contraction determined by the Stein factorization 
of $S'\to Z$. 
Write the pullback of $K_{S^\nu}+\Delta_{S^\nu}$ to $S'$ as $K_{S'}+\Delta_{S'}$. 
Since $-(K_X+\Delta)$ is ample$/Z$, we can see that  
$-(K_{S'}+\Delta_{S'})$ is nef and big$/V$. So running an LMMP$/V$ on $K_{S'}$ ends with a Mori fibre 
space over $V$ which is either $\PP^2$ or a $\PP^1$-bundle. So there is a covering family of rational curves 
 on $S'$ over $V$ such that for the general member $\Gamma_{S'}$ we have $-3 \le K_{S'}\cdot \Gamma_{S'}<0$, hence 
$-3 \le (K_{S'}+\Delta_{S'})\cdot \Gamma_{S'}<0$. 
Taking the image of the family on $S$ gives a covering family of curves on $S$ over $Z$ such that 
for a general member $\Gamma$ we have 
$$
-3\le (K_{S'}+\Delta_{S'})\cdot \Gamma_{S'}=(K_X+\Delta)\cdot \Gamma <(K_X+B)\cdot \Gamma <0
$$
so we are done in the divisorial case.

Now assume that $X\to Z$ is a flipping contraction and let $X\bir X^+/Z$ be its log flip 
which exists by [\ref{B}, Theorem 1.1]. Let $P^+$ be a sufficiently ample divisor on $X^+$.  
Let $\phi\colon W\to X$ be a log resolution of $(X,B)$ such that the induced map $\psi\colon W\bir X^+$ 
is a morphism. 
Let $B_W$ be the sum of the birational transform of $B$ and the reduced exceptional divisor of 
$\phi$, let $A_W$ be the pullback of $A$, and let $P_W$ be the pullback of $P^+$.
Since $(X,B)$ is klt, 
$$
K_W+B_W+A_W=\phi^*(K_X+B+A)+E_W\equiv E_W/Z
$$
where $E_W$ is effective and its support is equal to the support of the reduced exceptional divisor of $\phi$. 

Run an LMMP$/X$ on $K_W+B_W+A_W+P_W$. If $R_W$ is a 
$K_W+B_W+A_W+P_W$-negative extremal ray$/Z$, then $(K_W+B_W+A_W)\cdot R_W<0$ as $P_W$ is nef, 
hence $E_W\cdot R_W<0$. By [\ref{B}, Theorem 1.5], $R_W$ can be contracted and all the 
curves generating $R_W$ are inside $\Supp E_W$. 
Therefore, by restricting to the components of $\rddown{B_W}$ 
and applying the cone theorem on surfaces, we can find a  
curve $\Gamma_W$ generating $R_W$ such that 
$$
-3\le (K_W+B_W+A_W)\cdot \Gamma_W<0
$$  
In particular, $P_W\cdot \Gamma_W\le 3$.
Since $P_W$ is the pullback of a sufficiently ample divisor, it is necessary to have $P_W\cdot \Gamma_W=0$. 
Therefore, $R_W$ is contracted over $X^+$. 
Arguing in the same way, one shows that $P_W$ intersects every extremal ray in the 
process of the LMMP trivially. 
Moreover, by special termination [\ref{B}, Proposition 5.5], the LMMP terminates with a model $Y/X$ 
and the induced map $Y\bir X^+$ is a morphism. 
Denote $Y\to X$ and $Y\to X^+$ by $\alpha$ and $\beta$ respectively.

By construction, 
$$
K_Y+B_Y+A_Y\equiv E_Y/Z 
$$ 
where $E_Y$ is effective and its support is equal to the support of the reduced exceptional divisor of 
$\alpha$, which is equal to the reduced exceptional divisor of $\beta$. 
Obviously, $\beta$ is not an isomorphism, hence it contracts some divisor as $X^+$ is $\Q$-factorial. 
Thus $E_Y\neq 0$, and by the negativity lemma, $K_Y+B_Y+A_Y+P_Y$ is not nef$/X^+$. 
On the other hand, if $a\gg 0$, then $K_Y+B_Y+a A_Y+P_Y$ is nef$/Z$ because arguing as in the 
last paragraph we know that any $K_Y+B_Y+A_Y+P_Y$-negative extremal ray$/Z$ is generated 
by a  curve $C$ with 
$$
-3 \le (K_Y+B_Y+A_Y+P_Y)\cdot C
$$ 
By construction of $Y$, $C$ cannot be contracted over $X$ and so $A_Y\cdot C>0$.
So we have shown that 
 $K_Y+B_Y+A_Y$ is not nef$/X^+$ but $K_Y+B_Y+a A_Y$ is nef$/X^+$ for any $a\gg 0$.

Let $\lambda$ be the smallest number such that $K_Y+B_Y+\lambda A_Y$ is nef$/X^+$.
Note that $\lambda>1$ since $K_Y+B_Y+A_Y$ is not nef$/X^+$. 
Now by [\ref{B}, 3.4] there is an extremal ray $R_Y/X^+$ such that 
$(K_Y+B_Y)\cdot R_Y<0$ but $(K_Y+B_Y+\lambda A_Y)\cdot R_Y=0$. By construction, 
$E_Y\cdot R_Y<0$, so there is a rational curve $\Gamma_Y$ generating $R_Y$ such that 
$-3\le (K_Y+B_Y)\cdot \Gamma_Y$. Thus $A_Y\cdot \Gamma_Y\le 3$ as $(K_Y+B_Y+A_Y)\cdot \Gamma_Y<0$. 
Since $K_Y+B_Y+A_Y$ is nef$/X$, $\Gamma_Y$ is not contracted over $X$.
Let $\Gamma\subset X$ be the image of $\Gamma_Y$. Then $\Gamma$ generates $R$ and 
$A\cdot \Gamma\le A_Y\cdot \Gamma_Y\le 3$. Therefore, $-3\le (K_X+B)\cdot \Gamma$.\\   
\end{proof}

\begin{lem}\label{l-nef-thresh}
Let $(X,B)$ be a $\Q$-factorial projective dlt pair of dimension $3$ over 
$k$ of char $p>5$ such that $B$ is a $\Q$-boundary and $K_X+B$ is not nef. 
Then there is a natural number $n$ depending only on $(X,B)$ such that if $H$ is 
an ample Cartier divisor and 
$$
\lambda=\min\{t \mid K_X+B+tH~~\mbox{is nef}\}
$$   
then $\lambda=\frac{n}{m}$ for some natural number $m$, and there is a rational curve $\Gamma$
such that 
$$
-6\leq (K_X+B)\cdot\Gamma\leq 0 ~~\mbox{and}~~ (K_X+B+\lambda H)\cdot\Gamma= 0
$$
\end{lem}
\begin{proof}
Let $I$ be a natural number so that $I(K_X+B)$ is Cartier. Fix an ample Cartier divisor $H$ 
and let $\lambda$ be as in the statement of the proposition.
It is enough to show that there is a rational curve $\Gamma$ satisfying the last claim of the 
proposition because  then    
$\lambda=\frac{-I(K_X+B)\cdot \Gamma}{IH\cdot \Gamma}$ so taking $n=(6I)!$ we can write 
$\lambda=\frac{n}{m}$ for some natural number $m$.

Now assume that $L=K_X+B+\lambda H$ is big. Then by [\ref{B}, 3.4], 
there is an extremal ray $R$ such that $L\cdot R=0$. Moreover, by [\ref{B}, 3.3], 
$R$ is generated by some curve, hence $\lambda$ is rational and $L$ is $\Q$-Cartier. 
Since $L$ is nef and big, we can apply Lemma \ref{l-f-d-ray}. 

From now on we assume $L$ is not big. By extending $k$ we can assume it is uncountable. 
By Theorem \ref{t-nef-dimension}, the nef dimension of $L$ is at most $2$ and there is a nef reduction 
map $X\bir Z$ for $L$. For the moment assume that $\dim Z>0$.
Let $\phi\colon W\to X$ be a resolution so that $h\colon W\bir Z$ is a morphism. 
The map $X\bir Z$ is regular and proper over some open subset $V\subseteq Z$. 
Let $P$ be a general effective Cartier divisor on 
$Z$ intersecting $V$ and let $G$ be its pullback to $X$, that is, $G=\phi_*h^*P$. 
Let $S$ be a component of $G$ whose generic point maps into $V$ and such that $S$ is not a 
component of $\rddown{B}$, and let $Q$ be the image of $S$ on $Z$.  
There is $s>0$ such that the coefficient of $S$ in 
$B+sG$ is $1$. Let $\Theta=B+sG$ and let $S^\nu$ be the normalization of $S$ (note that 
$\Theta$ is not necessarily a boundary). By adjunction, we can write 
$$
K_{S^\nu}+\Theta_{S^\nu}=(K_X+\Theta)|_{S^\nu}
$$
for some $\Theta_{S^\nu}\ge 0$. Let $H_{S^\nu}=H|_{S^\nu}$. Then $K_{S^\nu}+\Theta_{S^\nu}+\lambda H_{S^\nu}$ 
is numerically trivial over the generic point of $V_Q:=V\cap Q$. 
Let $T\to {S^\nu}$ be the minimal resolution of $S^\nu$. Then the pullback of 
 $K_{S^\nu}+\Theta_{S^\nu}+\lambda H_{S^\nu}$ to $T$ can be written as $K_T+\Theta_T+\lambda H_T$ where $H_T$ 
 is the pullback of $H_{S^\nu}$. Run the LMMP$/V_Q$ on $K_T$ which ends with a Mori fibre space 
 because $K_T\equiv -\Theta_T-\lambda H_T$ over the generic point of $V_Q$ and because $H_T$ is big. The Mori fibre space is either 
 $\PP^2$ or a $\PP^1$-bundle.
Therefore, there is a covering family of rational curves on $T/Q$ such that $-3\le K_{T}\cdot \Gamma_{T}<0$ 
for the general members $\Gamma_{T}$ of the family. Since $\Theta_T\cdot \Gamma_T\ge 0$ and $H_T\cdot \Gamma_T>0$, 
we get $-3\le (K_T+\Theta_T)\cdot \Gamma_T<0$ for the very general members $\Gamma_T$. 
Taking the image of the family on 
$X$ we get a covering family of rational curves of $S/Q$ such that 
$$
-3\le (K_T+\Theta_T)\cdot \Gamma_T=(K_X+\Theta)\cdot \Gamma=(K_X+B)\cdot \Gamma<0 
$$  
for the very general members $\Gamma$ of the family. Note that by construction, $L\cdot \Gamma=0$.

Finally we treat the case $\dim Z=0$, that is, when $L\equiv 0$. 
In this case $-(K_X+B)$ is ample. Let $C$ be a smooth projective curve inside the smooth locus of 
$X$ such that $B\cdot C\ge 0$. We can obtain such $C$ by cutting $X$ by hypersurface sections. 
Note that $K_X\cdot C<0$.
 Fix a closed point $c\in C$. 
Now by [\ref{kollar}, Chapter II, Theorem 5.8], there is a rational curve $\Gamma$ passing 
through $c$ such that 
$$
\lambda H\cdot \Gamma \le 6\frac{\lambda H\cdot C}{-K_X\cdot C } 
$$  
Since $B\cdot C\ge 0$, we have $-K_X\cdot C\ge -(K_X+B)\cdot C$, hence 
$$
\lambda H\cdot \Gamma \le 6\frac{\lambda H\cdot C}{-(K_X+B)\cdot C }=6 
$$
But then $-(K_X+B)\cdot \Gamma=\lambda H\cdot \Gamma\le 6$ as required.\\
\end{proof}

\begin{proof}(of Theorem \ref{t-cone})
First we prove that the $K_X+B$-negative extremal rays do not accumulate in $\overline{NE}(X)_{K_X+B<0}$. 
Assume that there is a sequence $R_i$ of $K_X+B$-negative 
extremal rays which accumulate to some $K_X+B$-negative ray (not necessarily extremal).
Replacing the sequence and perturbing the coefficients of $B$ we can assume $B$ is a  $\Q$-boundary. 
By Lemma \ref{l-nef-thresh} and [\ref{KM}, Theorem 3.15], there is a collection of rays $\tilde{R}_j$ 
of $\overline{NE}(X)$ such that 
$$
\overline{NE}(X)=\overline{NE}(X)_{K_X+B\ge 0}+\sum_j \tilde{R}_j
$$ 
and such that the $\tilde{R}_j$ do not accumulate in $\overline{NE}(X)_{K_X+B<0}$.
For each $i$, since $R_i$ is extremal, there is some $j$ such that $R_i=\tilde{R}_j$. 
Therefore, the $R_i$ cannot accumulate in $\overline{NE}(X)_{K_X+B< 0}$, a contradiction.

Next we prove that there are only finitely many $K_X+B+A$-negative extremal rays for any 
ample $\R$-divisor $A$. 
Assume that there is an infinite sequence $R_i$ of $K_X+B+A$-negative extremal rays. 
Replacing the sequence we can assume the limit of $R_i$ exists as a ray, say $R$. 
By the last paragraph, $(K_X+B+A)\cdot R=0$. But then $(K_X+B)\cdot R<0$, hence the $R_i$ are an accumulating sequence of 
$K_X+B$-negative extremal rays which contradicts the last paragraph.

Now let $R$ be a $K_X+B$-negative extremal ray. We will show that $R$ is generated by a
 rational $\Gamma$ such that $-6\le (K_X+B)\cdot \Gamma$. Let $A$ be an ample $\R$-divisor. 
 Pick  $\epsilon,\delta>0$ so that 
$(K_X+B+(\epsilon+\delta) A)\cdot R<0$.  Since there are only finitely many $K_X+B+\epsilon A$-negative 
extremal rays, we can find a nef $\R$-divisor $N$ such that $N^\perp=R$. 
In particular, $R$ is the only $K_X+B+nN+\epsilon A$-negative extremal ray if $n$ is large 
enough. Now there is an ample $\R$-divisor 
$A'$ with sufficiently small coefficients and supported on $\Supp (B+nN+\epsilon A)$ 
so that there is a $\Q$-boundary 
$B'\sim_\Q B+nN+\epsilon A+A'$ with $(X,B')$ dlt and $(K_X+B'+\delta A)\cdot R<0$.   
Therefore, by Lemma \ref{l-nef-thresh}, we can find an ample $\Q$-divisor $H'$ so that 
$L'=K_X+B'+H'$ is nef and $L'^\perp=R$. Moreover, by Lemma \ref{l-f-d-ray} (if $L'$ is big) 
and  by Lemma \ref{l-nef-thresh} (if $L'$ is not big), 
there is a rational curve $\Gamma$ generating $R$ such that $-6\le (K_X+B')\cdot \Gamma$. 
From $(K_X+B'+\delta A)\cdot \Gamma<0$, we get $\delta A\cdot \Gamma\le 6$, hence 
$\Gamma$ belongs to a bounded family of curves, so there are only finitely many possibilities 
for $(K_X+B)\cdot \Gamma$. Moreover, by choosing $\epsilon$ and the coefficients of $A'$ 
to be small enough, we can assume $(A'+\epsilon A)\cdot \Gamma$ is sufficiently small. 
On the other hand, 
$-6-(A'+\epsilon A)\cdot \Gamma\le (K_X+B)\cdot \Gamma$. Therefore, $-6\le (K_X+B)\cdot \Gamma$. 
This completes the proof of the theorem.\\
\end{proof}

\subsection{Lifting curves birationally}
Here we prove some results which we will need in the next subsection. 

\begin{lem}\label{l-lift-curve-bir}
Let $(X,B)$ be a $\Q$-factorial dlt pair of dimension $3$ over $k$ of char $p>5$. 
Assume $f\colon X\to Z$ is a $K_X+B$-negative extremal birational contraction such that 
$-S$ is ample$/Z$ for some component $S$ of $\rddown{B}$. Let $C$ be a curve on $Z$. 
Then there is a curve $D$ on $X$  such that the induced morphism $D\to Z$ 
maps $D$ birationally onto $C$. 
\end{lem} 
\begin{proof}
If $C$ is not inside the image of the exceptional locus of $f$, eg a flipping contraction, 
then the statement is clear. So we can assume $S$ is contracted and mapped onto $C$. 
By [\ref{HX}][\ref{B}, Lemma 5.2], $S$ is normal.  Since $f$ has connected fibres and since all the positive dimensional  
fibres are inside $S$, the fibres of $S\to C$ are also connected. 
Let $S\to C'$ be the contraction given by the Stein factorization of $S\to C$. Then $C'\to C$ 
is the normalization of $C$ and it is birational. 

By adjunction, we can write $K_S+B_S=(K_X+B)|_S$ where $(S,B_S)$ is dlt. Moreover, 
$-(K_S+B_S)$ is ample$/C'$, hence $-K_S$ is big$/C'$. Let $U$ be an open subset of $S$ 
such that $U$ is smooth and $U\to C'$ is proper over its image, say $V$. 
Running an LMMP on $K_{U}$ ends with a $\PP^1$-bundle $T\to V$. In particular, there is a 
curve $D_T$ on $T$ which maps birationally onto $V$. Now let $D\subset S$ be the 
birational transform of $D_T$. Then $D$ maps birationally onto $C$.\\ 
\end{proof}

\begin{lem}\label{l-lift-curve-bir2}
Let $(X,B)$ be a klt pair of dimension $3$ over $k$ of char $p>5$, and $C$ a curve on $X$. 
Let $\phi\colon W\to X$ be a log resolution of $(X,B)$. Then there is a curve $D$ on $W$ such that 
the induced map $D\to X$ maps $D$ birationally onto $C$.
\end{lem}
\begin{proof}
Let $B_W$ be the sum of the birational transform of $B$ and the reduced exceptional 
divisor of $\phi$. Then $K_W+B_W=\phi^*(K_X+B)+E$ where $E$ is effective and 
its support is equal to the reduced exceptional divisor of $\phi$. Run an LMMP$/X$ 
on $K_W+B_W$ which is also an LMMP on $E$ whose support is inside $\rddown{B_W}$. 
By special termination [\ref{B}, Proposition 5.5], the LMMP 
ends with a model $Y/X$, and by the negativity lemma $E$ is contracted, hence 
$Y\to X$ is a small contraction. There is a curve $D_Y$ on $Y$ mapping birationally onto $C$. 
Let $W_i\bir W_{i+1}/Z_i$ be a step of the LMMP which is either a flip or a 
divisorial contraction with $W_{i+1}=Z_i$. 
Assume we have already found a curve $D_{i+1}$ on $W_{i+1}$ mapping birationally onto $C$. 
In the divisorial contraction case,  
apply Lemma \ref{l-lift-curve-bir} to find  a curve $D_{i}$ on $W_{i}$ mapping birationally onto $D_{i+1}$, 
hence mapping birationally onto $D$. In the flip case, apply Lemma \ref{l-lift-curve-bir} to find a 
curve $D_{i}$ on $W_{i}$ mapping birationally onto the image of $D_{i+1}$ on $Z_i$, so it 
also maps birationally onto $C$. So inductively we can find the required $D$ on $W$.\\  
\end{proof}

\subsection{Polytopes of boundary divisors}\label{ss-polytopes}
In this subsection, we study polytopes of divisors. 
This is similar to the characteristic $0$ case as treated in [\ref{B-mmodel-II}, Section 3] 
but for convenience we reproduce the details. 

We first recall the definition of extremal curves. 
Let $X$ be a normal projective variety and $H$ a fixed ample Cartier divisor 
(in practice we do not mention $H$ and assume that it is already fixed). 
Let $R$ be a ray of  $\overline{NE}(X)$. An \emph{extremal curve} for $R$ 
is a curve $\Gamma$ generating $R$ such that $H\cdot \Gamma\le H\cdot C$ 
for any other curve $C$ generating $R$. 
Let $D$ be an $\R$-Cartier divisor with $D\cdot R<0$. 
Then $D\cdot \Gamma\ge D\cdot C$ for any other curve $C$ generating $R$ [\ref{B}, 3.1].

Now let $X$ be a $\Q$-factorial projective klt variety of dimension $3$ over $k$ of char $p>5$. 
Let ${V}$ be a finite-dimensional rational affine subspace of the space of Weil $\R$-divisors on $X$. 
As mentioned in the introduction,
$$
\mathcal{L}=\{\Delta \in V \mid (X,\Delta) ~~\mbox{is lc}\}
$$
is a rational polytope in $V$. By  Theorem \ref{t-cone},  
for any $\Delta\in \mathcal{L}$ and any extremal curve $\Gamma$ of an extremal ray $R$ 
we have $(K_X+\Delta)\cdot \Gamma \ge -6$; note that although 
$(X,\Delta)$ may not be dlt but we can use the fact that $(X,a\Delta)$ is klt for any 
$a\in [0,1)$ and then take the limit over $a$.

Let $B_1,\dots, B_r$ be the vertices of $\mathcal{L}$, and let $m\in \N$ such that $m(K_X+B_j)$ are Cartier. 
For any $B\in \mathcal{L}$, there are  
 nonnegative real numbers $a_1,\dots, a_r$ such that $B=\sum a_jB_j$ and $\sum a_j=1$. 
Moreover, for any curve $\Gamma$ on $X$ 
the intersection number 
$$
(K_X+B)\cdot \Gamma=\sum a_j(K_X+B_j)\cdot \Gamma
$$ 
is of the form 
$\sum a_j\frac{n_j}{m}$ for certain $n_1,\dots, n_r\in \Z$. 
If $\Gamma$ is an extremal curve of an extremal ray, then the $n_j$ satisfy $n_j\ge -6m$.

For an $\R$-divisor $D=\sum d_iD_i$ where the $D_i$ are the irreducible components of $D$, we define $||D||:=\max\{|d_i|\}$.

\begin{prop}\label{p-polytope-rays}
Let $X$, $V$, and $\mathcal{L}$ be as above, and fix $B\in \mathcal{L}$. Then  
there are real numbers $\alpha,\delta>0$, depending only on $(X,B)$ and $V$, such that 
\begin{enumerate}
\item if $\Gamma$ is any extremal curve and if 
$(K_X+B)\cdot \Gamma>0$, then $(K_X+B)\cdot \Gamma>\alpha$;

\item if $\Delta\in \mathcal{L}$, $||\Delta-B||<\delta$ and  $(K_X+\Delta)\cdot R\le 0$ for an extremal ray $R$, 
then $(K_X+B)\cdot R\le 0$;

\item  let $\{R_t\}_{t\in T}$ be a family of extremal rays of $\overline{NE}(X)$. Then the set 
$$ 
\mathcal{N}_T=\{\Delta \in \mathcal{L} \mid (K_X+\Delta)\cdot R_t\ge 0 ~~\mbox{for any $t\in T$}\}
$$
is a rational polytope;

\item assume $K_X+B$ is nef, $\Delta\in \mathcal{L}$, and that $X_i\bir X_{i+1}/Z_i$ is a  
sequence of $K_X+\Delta$-flips which are $K_X+B$-trivial and $X=X_1$; then for any curve $\Gamma$ on any $X_i$, we have 
$(K_{X_i}+B_i)\cdot \Gamma>\alpha$ if $(K_{X_i}+B_i)\cdot \Gamma>0$ where $B_i$ is the 
birational transform of $B$;

\item in addition to the assumptions of $(4)$ suppose  that $||\Delta-B||<\delta$; if 
$(K_{X_i}+\Delta_i)\cdot R\le 0$ for an extremal ray $R$ on some $X_i$, then $(K_{X_i}+B_i)\cdot R= 0$ 
where $\Delta_i$ is the birational transform of $\Delta$. 
\end{enumerate}
\end{prop}
\begin{proof}
Let $B_1,\dots, B_r$ be the vertices of $\mathcal{L}$ and  
let $m$ be a natural number so that $m(K_X+B_j)$ are all Cartier. 
Write $B=\sum a_jB_j$ where $a_j\ge 0$ and $\sum a_j=1$.\\ 

(1) If $B$ is a $\Q$-divisor, then the statement is trivially true even if $\Gamma$ is not extremal. If $B$ 
is not a $\Q$-divisor, then  
$$
(K_X+B)\cdot \Gamma=\sum a_j(K_X+B_j)\cdot \Gamma
$$ 
and if $(K_X+B)\cdot \Gamma<1$, then there are only finitely many possibilities for the intersection 
numbers $(K_X+B_j)\cdot \Gamma$ because $(K_X+B_j)\cdot \Gamma \ge -6$ as $\Gamma$ is extremal, 
and this in turn implies that 
there are only finitely many possibilities for the intersection 
number $(K_X+B)\cdot \Gamma$. So 
the existence of $\alpha$ is clear.\\ 

(2)  If the statement is not 
true then there is an infinite sequence of $\Delta_t\in \mathcal{L}$ and extremal rays $R_t$ such that for 
each $t$ we have 
$$
(K_X+\Delta_t)\cdot R_t\le 0  ~~~~~~~~ \mbox{but}~~~~ ~~~~~ (K_X+B) \cdot R_t> 0
$$ 
and $||\Delta_t-B||$ converges to $0$.  
There are nonnegative real numbers $a_{1,t},\dots, a_{r,t}$ such that $\Delta_t=\sum a_{j,t} B_j$ and $\sum a_{j,t}=1$. 
Since $||\Delta_t-B||$ converges to $0$, $a_j=\lim_{t\to \infty} a_{j,t}$. Perhaps after replacing the sequence with an infinite subsequence 
we can assume that the sign of $(K_X+B_j)\cdot R_t$ is independent of $t$ for each $j$. 
On the other hand, we can assume that for each $t$ there is $\Delta_t'\in \mathcal{L}$ such that $(K_X+\Delta_t')\cdot R_t<0$, 
hence  we have an extremal 
curve $\Gamma_t$ for $R_t$, by Theorem \ref{t-cone}. 

Now, if $(K_X+B_j)\cdot \Gamma_t\le 0$, then $-6\le (K_X+B_j)\cdot \Gamma_t\le 0$ by Theorem \ref{t-cone}, 
hence there are only finitely many possibilities for $(K_X+B_j)\cdot \Gamma_t$, so  
we can assume that it is independent of $t$. On the other hand, if $(K_X+B_j)\cdot \Gamma_t> 0$ 
and if $a_j\neq 0$, then  $(K_X+B_j)\cdot \Gamma_t$ 
 is bounded from below and above because 
$$
(K_X+\Delta_t)\cdot \Gamma_t=\sum a_{j,t}(K_X+B_j)\cdot \Gamma_t \le 0
$$ 
and because for $t\gg 0$, $a_{j,t}$ is bounded from below as it is sufficiently close to $a_j$.  
Therefore, if $a_j\neq 0$, then there are only finitely many possibilities for  $(K_X+B_j)\cdot \Gamma_t$  and 
we could assume that it is independent of $t$. 

Rearranging the indexes we can assume that $a_j\neq 0$ for $1\le j\le l$ but $a_j=0$ for 
$j>l$. Then by (1) the number 
$$
(K_X+\Delta_t)\cdot \Gamma_t=\sum a_{j,t} (K_X+B_j)\cdot \Gamma_t =
$$
$$
 (K_X+B)\cdot \Gamma_t +\sum_{j\le l} (a_{j,t}-a_j) (K_X+B_j)\cdot \Gamma_t 
+\sum_{j>l} a_{j,t} (K_X+B_j)\cdot \Gamma_t
$$
is positive if $t\gg 0$ because $(K_X+B)\cdot \Gamma_t>\alpha$, and 
if $j\le l$, then $|(a_{j,t}-a_j) (K_X+B_j)\cdot \Gamma_t|$ is sufficiently 
small, and if $j>l$, then $a_{j,t} (K_X+B_j)\cdot \Gamma_t$ is either positive or 
$|a_{j,t} (K_X+B_j)\cdot \Gamma_t|$ is sufficiently 
small. This is a contradiction.\\

(3) We may assume that for each $t\in T$ there is some $\Delta\in \mathcal{L}$ such that $(K_X+\Delta)\cdot R_t<0$, 
in particular, $(K_X+B_j)\cdot R_t<0$ for some vertex $B_j$ of $\mathcal{L}$ and that 
$R_t$ is generated by some extremal curve. Since by Theorem \ref{t-cone} the set of such extremal rays is discrete, 
we may assume that $T\subseteq \N$.

Obviously, $\mathcal{N}_T=\bigcap_{t\in T} \mathcal{N}_{\{t\}}$ is a convex compact subset of 
$\mathcal{L}$ since each $\mathcal{N}_{\{t\}}$ is a convex closed subset. If $T$ is finite, the claim is trivial 
because $\mathcal{N}_T$ is then cut out of $\mathcal{L}$ by finitely many inequalities with rational coefficients. 
So we may assume that $T=\N$. By (2) and by the compactness of $\mathcal{N}_T$, there are 
$\Delta_1,\dots, \Delta_n\in \mathcal{N}_T$ 
and $\delta_1,\dots,\delta_n>0$ such that $\mathcal{N}_T$ is covered by the balls 
$\mathcal{B}_i=\{\Delta \in \mathcal{L} \mid ||\Delta-\Delta_i||<\delta_i\}$ and such that 
if $\Delta\in \mathcal{B}_i$ with $(K_X+\Delta)\cdot R_t<0$ for some $t$, then $(K_X+\Delta_i)\cdot R_t=0$.

Let
$$
T_i=\{t\in T \mid (K_X+\Delta)\cdot R_t<0 ~~\mbox{for some $\Delta \in \mathcal{B}_i$}\}
$$
Then by construction $(K_X+\Delta_i)\cdot R_t=0$ for any $t\in T_i$. 
We claim that  
$$
\mathcal{N}_T=\bigcap_{1\le i\le n} \mathcal{N}_{T_i}
$$
Let $T'=\bigcup_{1\le i\le n}  T_i$ and let $S=T\setminus T'$. Pick $s\in S$. Since $s\notin T_i$ for each $i$, 
$(K_X+\Delta)\cdot R_s\ge 0$ for every $\Delta \in \bigcup_{1\le i\le n} \mathcal{B}_i$. Thus 
$ \mathcal{N}_T\subseteq \bigcup_{1\le i\le n} \mathcal{B}_i\subseteq \mathcal{N}_S$ which in turn implies that  
$$
\mathcal{N}_{T}=\mathcal{N}_{T'}=\bigcap_{1\le i\le n} \mathcal{N}_{T_i}
$$
because the $\mathcal{B}_i$ give an open cover of $\mathcal{N}_T$ 
and $\mathcal{N}_{T'}$ is a convex closed set containing $\mathcal{N}_T$.

By the last paragraph, it is enough to prove that each $\mathcal{N}_{T_i}$ is a rational polytope and by replacing $T$ with $T_i$, 
we could assume from the beginning that there is some $\Delta\in \mathcal{N}_T$ such that 
$(K_X+\Delta)\cdot R_t=0$ for every $t\in T$.
 If $\dim \mathcal{L}=1$, this already proves the proposition. 

Assume $\dim \mathcal{L}>1$ and let $\mathcal{L}^1, \dots, \mathcal{L}^p$ be the proper faces of $\mathcal{L}$.
Then each $\mathcal{N}_{T}^i:=\mathcal{N}_{T}\cap \mathcal{L}^i$ is a rational polytope by induction. 
Moreover,  for each $\Delta''\in \mathcal{N}_{T}$ 
which is not equal to $\Delta$, there is $\Delta'$ on some face $\mathcal{L}^i$ such that $\Delta''$ is on the line 
segment determined by $\Delta$ and $\Delta'$. Since $(K_X+\Delta)\cdot R_t=0$ for every 
$t\in T$,  $\Delta'\in \mathcal{N}_T^i$.
Hence $\mathcal{N}_T$ is the convex hull of $\Delta$ and all the  $\mathcal{N}_{T}^i$. 
Now, there is a finite subset $V\subset T$ such that for each $i$ we have 
$\mathcal{N}_T^i=\mathcal{N}_{V}\bigcap \mathcal{L}^i$.
But then the convex hull of  $\Delta$ and all the $\mathcal{N}_{T}^i$ 
is nothing but $\mathcal{N}_{V}$, hence $\mathcal{N}_T=\mathcal{N}_V$ and we are done.\\

(4) Note that the sequence being $K_X+B$-trivial means that $K_{X_i}+B_i$ is numerically 
trivial over $Z_i$ for each $i$. This in particular implies that $K_{X_i}+B_i$ is nef for every $i$.
Since $K_X+B$ is nef, $B\in \mathcal{N}_T$ where we take $\{R_t\}_{t\in T}$ to be the family of 
all the extremal rays of $\overline{NE}(X)$. Since $\mathcal{N}_T$ is a rational polytope by (3), there are 
positive real numbers $a_1',\dots,a_{r'}'$, and $m'\in \N$ so that $\sum a_j'=1$, $B=\sum a_j'B_j'$, 
and each $m'(K_X+B_j')$ is Cartier  
where $B_j'$ are among the vertices of $\mathcal{N}_T$. Therefore, since 
$K_X+B=\sum a_j'(K_X+B_j')$ and since each $K_X+B_j'$ is nef, the sequence $X_i\bir X_{i+1}/Z_i$ 
is also $K_X+B_j'$-trivial for each $j$. Thus $K_{X_i}+B_{j,i}'$ is nef and $({X_i},B_{j,i}')$ 
is log canonical for every $i$. 

Now fix $i$ and let $\phi\colon W\to X$ and $\psi\colon W\to X_i$ 
be a common log resolution. Since $X=X_1$ is klt, $(X_i,\Theta)$ is also klt 
for some $\Theta$. 
Then by Lemma \ref{l-lift-curve-bir2}, there is a curve $D$ on $W$ which maps birationally 
onto $\Gamma\subset X_i$. This implies that if $(K_{X_i}+B_{j,i}')\cdot \Gamma >0$, then 
$$
(K_{X_i}+B_{i})\cdot \Gamma\ge a_j'(K_{X_i}+B_{j,i}')\cdot \Gamma=a_j'\psi^*(K_{X_i}+B_{j,i}')\cdot D= a_j'\phi^*(K_X+B_j')\cdot D\ge \frac{a_j'}{m'}
$$ 
 Therefore, perhaps after replacing $\alpha$ of $(1)$ 
with a smaller one, we have $(K_{X_i}+B_{i})\cdot \Gamma >\alpha$ if $(K_{X_i}+B_{i})\cdot \Gamma >0$.\\

(5) Let $\Delta'$ be on the boundary of $\mathcal{L}$ so that $\Delta$ belongs to the line segment 
determined by $B$ and $\Delta'$.  There are nonnegative real numbers 
$r,s$ such that $s>0$, $r+s=1$ and $\Delta=rB+s\Delta'$.  In particular, the sequence $X_i\bir X_{i+1}/Z_i$ is also a 
sequence of $K_X+\Delta'$-flips and $(X_i,\Delta_i')$ is lc for each $i$. 
Suppose that there is an extremal ray $R$ on some $X_i$ such that $(K_{X_i}+\Delta_i)\cdot R\le 0$ 
but $(K_{X_i}+B_i)\cdot R>0$.  By Theorem \ref{t-cone}, $(K_{X_i}+\Delta_i')\cdot \Gamma \ge -6$ 
for some curve $\Gamma$ generating $R$ (note that $X_i$ is $\Q$-factorial klt so we can 
apply \ref{t-cone} to $(X_i,\Delta_i')$ although it may not be dlt) On the other hand, by (4), $(K_{X_i}+B_i)\cdot \Gamma >\alpha$. Now 
$$
(K_{X_i}+\Delta_i)\cdot \Gamma=r(K_{X_i}+B_i)\cdot \Gamma+s(K_{X_i}+\Delta_i')\cdot \Gamma>r\alpha-6s
$$
and it is obvious that this is positive if $r>\frac{6s}{\alpha}$. In other words, if $\Delta$ is sufficiently close to 
$B$, then we get a contradiction. Therefore, it is enough to replace the $\delta$ of (2) by one sufficiently smaller.\\ 
\end{proof}

\subsection{Big log divisors}
We can derive base point freeness in the big case from results we have proved so far.

\begin{prop}\label{p-bpf-big}
Let $(X,B)$ be a projective klt pair of dimension $3$ and $X\to Z$ a projective contraction 
over $k$ of char $p>5$. 
If $K_X+B$ is nef and big$/Z$, then it is semi-ample$/Z$.
\end{prop}
\begin{proof}
Let $P$ be the pullback of a sufficiently ample divisor on $Z$. By Theorem \ref{t-cone}, 
$K_X+B+P$ is globally nef and big. Since $P$ is nef and $K_X+B+P$ is nef and big, 
there exist $\epsilon >0$ and 
$$
\Delta\sim_\R B+P+\epsilon(K_X+B+P)
$$  
such that $(X,\Delta)$ is klt. It is enough to show $K_X+\Delta$ is semi-ample. 
Replacing $B$ with $\Delta$ we can then assume $Z$ is a point. 

By Proposition \ref{p-polytope-rays} (3), there are $\Q$-boundaries $B_j$ and nonnegative 
real numbers $a_j$ such that $\sum a_j=1$, $K_X+B=\sum a_j(K_X+B_j)$, and such that $K_X+B_j$ is nef for each $j$. 
Since $(X,B)$ is klt and $K_X+B$ is big, we can choose the $B_j$ so that 
$(X,B_j)$ is klt and $K_X+B_j$ is big for each $j$. By [\ref{B}, Theorem 1.4], 
each $K_X+B_j$ is semi-ample, hence $K_X+B$ is also semi-ample.\\
\end{proof}

\section{Finiteness of minimal models and termination}\label{s-finiteness}

In this section we prove finiteness of minimal models and 
derive termination with scaling under certain assumptions.

\begin{rem}\label{r-local}
With the setting as in Theorem \ref{t-finiteness},  
let $B\in \mathcal{L}_{A}(V)$ such that $(X,B)$ is klt. 
We can write $A\sim_\R A'+G/Z$ where $A'\ge 0$ is an ample $\Q$-divisor 
and $G\ge 0$ is also a $\Q$-divisor. Then there is a sufficiently small rational number $\epsilon>0$  
such that 
$$
(X,\Delta_B:=B-\epsilon A+\epsilon A'+\epsilon G)
$$
 is klt. 
Note that 
$$
K_X+\Delta_B\sim_\R K_X+B/Z
$$
Moreover, there is an open neighborhood of $B$ in $\mathcal{L}_{A}(V)$ such that for any $B'$ in that neighborhood 
$$
(X,\Delta_{B'}:=B'-\epsilon A+\epsilon A'+\epsilon G)
$$ 
is klt too. In particular, if 
$\mathcal{C}\subseteq \mathcal{L}_{A}(V)$ is a rational polytope containing $B$, 
then perhaps after shrinking $\mathcal{C}$ (but preserving its dimension) we can assume 
that $\mathcal{D}:=\{\Delta_{B'} \mid B'\in \mathcal{C}\}$ is a rational 
polytope of klt boundaries in $\mathcal{L}_{\epsilon A'}(W)$ where $W$ is the rational 
affine space $V+(1-\epsilon)A+\epsilon G$.  
The point is that we can change $A$ 
and get an ample part $\epsilon A'$ in the boundary. So, when we are concerned with a problem locally around $B$ we feel free to assume that $A$ is actually ample by replacing it with $\epsilon A'$.
\end{rem}

\begin{prop}\label{p-finiteness}
Theorem \ref{t-finiteness} holds either 

$(1)$ if $K_X+\Delta$ is big over $Z$ for every $\Delta\in \mathcal{C}$, or 

$(2)$ if Theorem \ref{t-bpf} holds.
\end{prop}
\begin{proof}
 We may assume that the dimension of $\mathcal{C}$ is positive. We can also assume that 
the proposition holds for polytopes with smaller dimension.
Since $\mathcal{C}$ is compact, it is enough to prove the statement locally near a fixed 
$B\in \mathcal{C}$. If $K_X+B$ is not pseudo-effective$/Z$, then the same 
holds in a neighborhood of $B$ inside $\mathcal{C}$. So 
we may assume that $K_X+B$ is pseudo-effective$/Z$. Then $(X,B)$ 
has a log minimal model $(Y,B_Y)$ over $Z$ by Theorem \ref{t-mmodel}. 
Moreover, the polytope $\mathcal{C}$ determines a rational polytope $\mathcal{C}_Y$ of 
$\R$-divisors on $Y$ by taking birational transforms of elements of $\mathcal{C}$. If we shrink $\mathcal{C}$ 
around $B$ we can assume that for every $\Delta\in \mathcal{C}$ the log discrepancies satisfy 
$$
a(D,X,\Delta)<a(D,Y,\Delta_Y)
$$
for any prime divisor $D$ on $X$ which is contracted by $X\bir Y$. 
So for each $\Delta\in \mathcal{C}$, a log minimal model of 
 $(Y,\Delta_Y)$ over $Z$ is also a log minimal model of $(X,\Delta)$ over $Z$. 
Therefore, we can replace $(X,B)$ by $(Y,B_Y)$ and replace $\mathcal{C}$ by $\mathcal{C}_Y$, hence 
from now on assume that $K_X+B$ is nef$/Z$. Then $K_X+B$ is semi-ample$/Z$: in case $(1)$ we use 
Proposition \ref{p-bpf-big} and in case $(2)$ we use Theorem \ref{t-bpf}.
So $K_X+B$ defines a contraction 
$f\colon X\to S/Z$.  

 Now by our assumption at the beginning of this proof, we may assume that 
there are finitely many birational maps $\psi_j\colon X\bir Y_j/S$ such that 
for any $\Delta'$ on the boundary of $\mathcal{C}$ with $K_X+\Delta'$ pseudo-effective$/S$, there is $j$ such that 
$(Y_j,\Delta_{Y_j}')$ is a log minimal model of $(X,\Delta')$ over $S$.
On the other hand, by Proposition \ref{p-polytope-rays}, 
there is a sufficiently small $\epsilon>0$ such that for any $\Delta\in \mathcal{C}$ 
with $||B-\Delta||<\epsilon$, and any $j$, and any $K_{Y_j}+\Delta_{Y_j}$-negative extremal ray $R/Z$ we have the equality $(K_{Y_j}+B_{Y_j})\cdot R=0$. 
Note that all the pairs $(Y_j,B_{Y_j})$ are klt and $K_{Y_j}+B_{Y_j}\equiv 0/S$ and nef$/Z$.

Pick $\Delta\in \mathcal{C}$ with $0<||B-\Delta||<\epsilon$ such that $K_X+\Delta$ is pseudo-effective$/Z$, and let $\Delta'$ be the unique point on 
the boundary of $\mathcal{C}$ such that $\Delta$ belongs to the line segment given by $B$ and $\Delta'$. Since $K_X+B\equiv 0/S$, there is some $t>0$ such that 
$$
K_X+\Delta'=K_X+B+\Delta'-B\equiv \Delta'-B=t(\Delta-B)\equiv t(K_X+\Delta)/S 
$$
hence 
$K_X+\Delta'$ is pseudo-effective$/S$, and $(Y_j,\Delta_{Y_j}')$ is a log minimal model of $(X,\Delta')$ over $S$ for some $j$. Moreover, $(Y_j,\Delta_{Y_j})$ is a log minimal model of $(X,\Delta)$ over $S$ for the same $j$. Furthermore, $(Y_j,\Delta_{Y_j})$ is a log minimal model of $(X,\Delta)$ over $Z$ because any $K_{Y_j}+\Delta_{Y_j}$-negative extremal ray $R/Z$ would be over $S$ by the last paragraph. Finally, we just need to shrink $\mathcal{C}$ around $B$ appropriately.\\
\end{proof}

\begin{prop}\label{p-term}
Theorem \ref{t-term} holds either 

$(1)$ if $K_X+B$ is big or 

$(2)$ if Theorems \ref{t-bpf} and \ref{t-contraction} hold. 
\end{prop}
\begin{proof}
Since $B$ is big$/Z$, we can assume $B\ge A$ for some (globally) ample $\R$-divisor $A$.
Thus there are only finitely many $K_X+B$-negative extremal 
rays by Theorem \ref{t-cone} (iii). Moreover, by adding to $A$ the pullback of a sufficiently 
ample divisor on $Z$ and changing $A$ up to $\R$-linear equivalence and applying Theorem 
\ref{t-cone}, we can assume that all the $K_X+B$-negative extremal rays are over $Z$ and that 
$K_X+B+C$ is globally nef. Therefore, 
we can run the LMMP globally which would automatically be over $Z$. 

Let $\lambda\ge 0$ be the smallest number such that $K_X+B+\lambda C$ is nef. We can assume $\lambda>0$.  
 So there is an extremal ray $R$ 
such that $(K_X+B)\cdot R<0$ but $(K_X+B+\lambda C)\cdot R=0$.
The ray $R$ can be contracted: in case $(1)$ we can use [\ref{B}, Theorem 1.5] but in case $(2)$ 
we use Theorem \ref{t-contraction}. If this contraction is not birational, 
then we have a Mori fibre space and we stop. 
Otherwise, we let $X\bir X'$ be the corresponding flip or divisorial contraction, and we continue 
with $X'$ and so on. This shows that we can run the LMMP with scaling. 

We will show that the LMMP terminates. Assume not. Assume that we get an infinite sequence 
$X_i\bir X_{i+1}/Z_i$ of log flips. We may assume that $X=X_1$. 
Let $\lambda_i$ be the numbers appearing in the LMMP with scaling, and put $\lambda=\lim \lambda_i$. 
By definition, $K_{X_i}+B_i+\lambda_iC_i$ is nef but numerically zero over $Z_i$ 
where $B_i$ and $C_i$ are the birational transforms of $B$ and $C$ respectively. 
Replacing $B$ with $B+\lambda C$, we can assume $\lambda=0$.

Let $H_1, \dots, H_m$ be general effective ample Cartier divisors 
on $X$ which generate the space $N^1(X)$. Since $B$ is big, we may assume that 
$B-\epsilon (H_1+\cdots+H_m)\ge 0$ for some rational number $\epsilon>0$. Replacing $A$ we 
can assume $A=\frac{\epsilon}{2} (H_1+\cdots+H_m)$.
Let $V$  be the $\R$-vector space generated by  
the components of $B+C$, and let $\mathcal{C}\subset \mathcal{L}_A(V)$ be a rational polytope of maximal dimension containing 
an open neighborhood of $B$ and such that $(X,\Delta)$ is klt for every $\Delta\in \mathcal{C}$. 
In case $(1)$ we can choose $\mathcal{C}$ so that  $K_X+\Delta$ is big for every $\Delta\in \mathcal{C}$. 

We can choose $\mathcal{C}$ such that for each $i$ there is an ample $\Q$-divisor $G_i=\sum g_{i,j}H_{i,j}$ on $X_i$ 
with sufficiently small coefficients, 
where $H_{i,j}$ on $X_i$ is the birational transform of $H_j$, such that $\Delta^i$ 
the birational transform of $B_i+G_i+\lambda_iC_i$ on $X$ belongs to $\mathcal{C}$. In particular, 
$K_{X_i}+B_i+G_i+\lambda_iC_i$ is ample and $(X_i,B_i+G_i+\lambda_iC_i)$ is the lc model of $(X,\Delta^i)$.

Now, by Proposition \ref{p-finiteness}, there are finitely many birational 
maps $\phi_l\colon X\bir Y_l$ such that for any $\Delta\in \mathcal{C}$ with $K_X+\Delta$ pseudo-effective, there 
is $l$ such that $(Y_l,\Delta_{Y_l})$ is a log minimal model of $(X,\Delta)$.
Since $K_{X_i}+B_i+G_i+\lambda_iC_i$ is ample and since the lc model is unique,
for each $i$, there is some $l$ such that $\phi_{1,i}\phi_l^{-1}$ is an isomorphism 
where $\phi_{i,j}$ is the birational map $X_i\bir X_j$. 
Therefore, there exist $l$ and an infinite set $I\subseteq \N$ such that $\phi_{1,i}\phi_l^{-1}$ is 
an isomorphism for any $i\in I$. This in turn implies that $\phi_{i,j}$ is an isomorphism for any $i,j\in I$. This is not possible as any log flip increases some log discrepancy.\\ 
\end{proof}

\begin{prop}\label{p-term-2}
Theorem \ref{t-term-2} holds if $K_X+B$ is pseudo-effective$/Z$.
\end{prop}
\begin{proof}
As in the proof of Proposition \ref{p-term}, we can reduce the problem to the global case, hence 
ignore $Z$. The fact that we can run the LMMP is also proved there. 
Let $X_i\bir X_{i+1}/Z_i$ be the steps of the LMMP which are either flips or divisorial contractions 
with $X_{i+1}=Z_i$, and $X=X_1$. Let $\lambda_i$ be the 
numbers that appear in the LMMP. Let $\lambda=\lim \lambda_i$. Assume that the LMMP 
does not terminate. We will derive a contradiction. 
If $\lambda>0$, then 
the LMMP is also an LMMP on $K_X+B+\frac{\lambda}{2}C$ with scaling of $(1-\frac{\lambda}{2})C$, 
so the theorem follows from Proposition \ref{p-term} in this case. Thus we can assume $\lambda=0$. 
In particular, this means $K_{X_i}+B_i$ is (numerically) a limit of movable divisors for 
any $i\gg 0$. Moreover, since $C$ is ample, we can assume that its components generate 
$N^1(X)$.

By Theorem \ref{t-mmodel}, $(X,B)$ has a log minimal model $(Y,B_Y)$.
Since $K_{X_i}+B_i$ is (numerically) a limit of movable divisors for 
any $i\gg 0$, the induced maps $X_i\bir Y$ are all isomorphisms in codimension one when $i\gg 0$. 
Fix $i\gg 0$. Since the components of $C_i$ generate 
$N^1(X_i)$, there is an ample divisor $H_i$ on $X_i$ supported on $\Supp C_i$.   
Let $H_Y$ be the birational transform of $H_i$ on $Y$. Pick a sufficiently small number $\epsilon>0$. 
Then  $(Y,B_Y+\epsilon H_Y +\lambda_i C_Y)$ is klt. Run an LMMP 
on $K_Y+B_Y+\epsilon H_Y +\lambda_i C_Y$ with scaling of some ample divisor. 
Since $\epsilon$ and $\lambda_i$ are sufficiently small and $K_Y+B_Y$ is nef, by Proposition \ref{p-polytope-rays} (5), 
the LMMP is $K_Y+B_Y$-trivial (note that although $H_Y$ depends on $i$ but $\Supp (\epsilon H_Y +\lambda_i C_Y)$ 
is independent of $i$ and its coefficients are sufficiently small). Moreover, the LMMP terminates by Proposition \ref{p-term}. 
As $K_{X_i}+B_i+\epsilon H_i+\lambda_i C_i$ is ample, the LMMP ends with $X_i$. 
Therefore, $K_{X_i}+B_i$ is also nef. This is a contradiction since $X_i\bir X_{i+1}/Z_i$ is a 
$K_{X_i}+B_i$-flip.\\ 
\end{proof}

\section{Relatively numerically trivial divisors}\label{s-num-trivial}

\begin{lem}\label{l-hyperplane-intersections}
Suppose $X$ is a normal projective $\Q$-factorial threefold over $k$  and $H_1,\dots,H_n$ are 
 ample divisors generating $N^1(X)$. Then the set of $1$-cycles $\{H_i\cdot H_j\}_{1\le i,j\le n}$ generates $N_1(X)$.
\end{lem}
\begin{proof}
Let ${V}$ be the vector subspace of $N_1(X)$ generated by $\{H_i\cdot H_j\}_{1\le i,j\le n}$.
Let $L$ be an $\R$-divisor with trivial intersection on every element of ${V}$. 
As the $H_i$ generate $N^1(X)$ we can write $L\equiv \sum_1^na_iH_i$.
Pick a curve $C$ on $X$ and let $S$ be a prime divisor on $X$ containing $C$.  Let ${S}^\nu$ be the normalization of $S$.  
We can write $S\equiv \sum_1^n b_iH_i$. Now 
 $$
 L|_S\cdot L|_S = L\cdot L\cdot S=L\cdot (\sum a_iH_i)\cdot (\sum b_iH_i) = 0
 $$  Similarly $L|_S\cdot H_i|_S=0$.  Take a resolution $\phi\colon S'\to {S}^\nu$. Apply the Hodge index theorem to
  $\phi^*(L|_{S^\nu})$: we can find an ample divisor $H'$ on $S'$ with $\phi_*H'=H_i|_{S^\nu}$ for some $i$, so  
$$
\phi^*(L|_{S^\nu})\cdot H'=\phi^*(L|_{S^\nu})\cdot \phi^*(H_i|_{S^\nu})=L|_S\cdot H_i|_S=0
$$ 
and also 
$\phi^*(L|_{S^\nu})\cdot\phi^*(L|_{S^\nu})=L|_S\cdot L|_S=0$.
Therefore, $\phi^*(L|_{S^\nu})\equiv 0$, hence $L|_S\num 0$ and $L\cdot C=0$. This shows that $L\equiv 0$. Since $L$ was chosen arbitrarily, 
$V=N_1(X)$.\\
\end{proof}

\begin{lem}\label{l-num-pullback}
Let $f\colon X\to Z$ be a projective 
contraction from a normal projective $\Q$-factorial threefold to a curve over $k$, and let $L$ be a nef $\Q$-divisor on $X$.
If $L\num 0/Z$, then $L\num f^*D$ for some $\Q$-divisor $D$ on $Z$.
\end{lem}
\begin{proof}
 Let $P$ be a point on $Z$.
Let $H_1,\dots, H_n$ be very ample divisors generating $N^1(X)$. 
We may assume each $H_i$ is normal and irreducible. By Lemma \ref{l-hyperplane-intersections}, 
the cycles $H_i\cdot H_j$ generate $N_1(X)$.    For each $i$, choose an ample divisor $A_i$ on $H_i$ and let $n_i$ be the value for which the divisor $D_i=n_iP$ on $Z$ satisfies $(L|_{H_i}-f|_{H_i}^*D_i)\cdot A_i=0$. 
Then  
$$
(L|_{H_i}-f|_{H_i}^*D_i)\cdot(L|_{H_i}-f|_{H_i}^*D_i) \geq 0
$$ 
because $L$ is nef, $(f|_{H_i}^*D_i)\cdot L|_{H_i} = 0$, and $(f|_{H_i}^*D_i)\cdot (f|_{H_i}^*D_i) = 0$.  Therefore by the Hodge index theorem, $L|_{H_i}-f|_{H_i}^*D_i\num 0$.  In particular, for each $i,j$ we have 
$$
(L-f^*D_i)\cdot H_i\cdot H_j = 0
$$  
hence 
$$
f^*D_i\cdot H_i\cdot H_j =f^*D_j\cdot H_i\cdot H_j 
$$ 
which implies that  $n_i=n_j$, so $D_i=D_j$. As the pairwise intersections $H_i\cdot H_j$ generate $N_1(X)$, 
we have $L\num f^*D$ where $D=D_i$.\\   
\end{proof}

\begin{lem}\label{l-relatively-trivial}
Let $f\colon X\to Z$ be a projective contraction from a normal quasi-projective variety onto a 
smooth curve over $k$. Assume $L$ is a nef$/Z$ $\R$-divisor on 
$X$ such that $L|_F\equiv 0$ where $F$ is the generic
fibre of $f$. Then $L\equiv 0/Z$.
\end{lem}
\begin{proof}
Pick a curve $C$ contracted by $f$. We will show that $L\cdot C=0$. 
Choose a surface $S$ containing $C$ such that $S\to Z$ is surjective. 
It is enough to show that $L|_S\cdot C=0$. So by replacing $X$ with the normalization $S^\nu$ 
and $X\to Z$ with the Stein factorization of $S^\nu\to Z$, we can assume 
$\dim X=2$. 

By extending $k$ we can assume it is uncountable. We will use some of the notation 
and results of [\ref{B-augmented}]. If $G$ is a very general fibre of $f$, then 
$h^0(\langle mL \rangle|_F)=h^0(\langle mL \rangle|_G)$ for every natural number $m$ 
where $F$ is the generic fibre. Since $L|_F$ is numerically trivial, 
$$
\limsup_{m\to +\infty} \frac{h^0(\langle mL \rangle|_F)}{m}=0
$$ 
by [\ref{B-augmented}, Proposition 4.3]. 
Thus 
$$
\limsup_{m\to +\infty} \frac{h^0(\langle mL \rangle|_G)}{m}=0
$$
hence $L|_G\equiv 0$ referring to the same result. 

Now let $G$ be a very general fibre over a closed point and $H$ any fibre over a closed point. 
Since $G\sim H$, we have $L\cdot H=L\cdot G=0$ which means $L|_H\equiv 0$. Therefore, 
$L\equiv 0/Z$ as claimed.\\  
\end{proof}

\begin{lem}\label{l-relatively-trivial-2}
Let $f\colon X\to Z$ be a flat projective morphism between quasi-projective schemes over $k$. 
Assume $L$ is a nef$/Z$ $\R$-divisor on $X$ and that $L|_F\equiv 0$ for the fibres $F$ of $f$ 
over some dense open subset of $Z$. Then $L\equiv 0/Z$. 
\end{lem}
\begin{proof}
Note that by saying $L$ is a nef$/Z$ $\R$-divisor on $X$ we mean $L$ is an element of 
$\Div(X)\otimes \R$ where $\Div(X)$ is the group of Cartier divisors on $X$ and 
$L\cdot C\ge 0$ for any curve $C$ contracted by $f$.  
  
Let $C$ be a curve contracted by $f$ to a point $z$. We want to show $L\cdot C=0$. 
We can assume there is a component $T$ of $Z$ of positive dimension containing $z$. 
Let $V$ be the normalization of a curve in $T$ passing through $z$. Let $Y=V\times_ZX$ and 
$g\colon Y\to V$ the induced morphism. Then $g$ is flat and we can choose $V$ so that 
$L|_Y$ is numerically trivial over some nonempty open subset of $V$. 
Replacing $f$ with $g$ we can 
assume $Z$ is a smooth  curve. Now the flatness means that every associated point of $X$  
maps to the generic point of $Z$, hence every irreducible component of $X$ 
maps onto $Z$. Replacing $X$ with one of its irreducible components 
(with reduced structure) containing $C$, we can assume 
$X$ is integral, that is, it is a variety. Now apply Lemma \ref{l-relatively-trivial} to the pullback of $L$ to 
the normalization of $X$.\\ 
\end{proof}

\begin{defn}
Suppose $f:X\to Z$ is a dominant morphism of normal varieties over $k$.  
 A divisor $L$ on $X$ is \emph{exceptional} if every component $L_i$ of $L$ is contracted by $f$, that is, 
$\dim f(L_i)<\dim(L_i)$, and $f(L_i)\neq Z$.
On the other hand, $L$ is \emph{very exceptional} if it is exceptional and for any prime divisor $P$ on $Z$, some divisorial component $Q$ of $f^{-1}(P)$ with $f(Q)=P$ is not contained in $\Supp L$.
\end{defn}

The following is an adaptation of  a result of Kawamata [\ref{Kawamata}, Proposition 2.1].

\begin{lem}\label{l-linear-pullback}
Let $f\colon X\to Z$ be a projection contraction between normal quasi-projective varieties over $k$
 and $L$ a nef$/Z$ $\R$-divisor on $X$ such that $L|_F\sim_\R 0$ where $F$ is the generic
fibre of $f$. Assume $\dim Z\le 3$ if $k$ has char $p>0$. Then there exist a diagram
$$
\xymatrix{
X'\ar[r]^\phi\ar[d]^{f'} & X\ar[d]^f\\
Z'\ar[r]^\psi & Z
}
$$
with $\phi,\psi$ projective birational, and an $\R$-Cartier divisor $D$ on $Z'$ such that 
$\phi^* L\sim_\R f'^*D$. 
Moreover, if $Z$ is $\Q$-factorial, then we can take $X'=X$ and $Z'=Z$. 
\end{lem}
\begin{proof}
Since $L|_F\sim_\R 0$, we can replace $L$ up to $\R$-linear equivalence hence assume 
 $L|_F=0$. So $L|_G=0$ for the general fibres $G$ of $f$. 
On the other hand, by flattening, there exist a diagram
$$
\xymatrix{
X''\ar[r]^\pi\ar[d]^{f''} & X\ar[d]^f\\
Z''\ar[r]^\mu & Z
}
$$
with $\pi,\mu$ projective birational and $f''$ flat (but $X''$ and $Z''$ may not be normal).
Applying Lemma \ref{l-relatively-trivial-2} to $\pi^*L$, we deduce that $\pi^*L\equiv 0/Z''$. 
Now we can extend the diagram as 
$$
\xymatrix{
X'\ar[r]\ar[d]^{f'} & X''\ar[r]\ar[d]^{f''} & X\ar[d]^f\\
Z'\ar[r] & Z''\ar[r] & Z
}
$$
with $X'\to X''$ and $Z'\to Z''$ projective birational, $X'$ normal and $Z'$ smooth. 
Denote the induced maps  $X'\to X$ and $Z'\to Z$ by $\phi$ and $\psi$ respectively.
Then $\phi^*L\equiv 0/Z'$. 

Since $\phi^*L$ is vertical$/Z'$,   
 there is an $\R$-divisor $E$ on $X'$ vertical$/Z'$ such that $\phi^*L+E\sim_\R 0/Z'$. 
For each divisorial component $P$ of $f(E)$, we can add $af^*P$ to $E$ for some 
appropriate number $a$ so that every component of $E$ mapping to $P$ has nonnegative 
coefficient, and at least one of them has coefficient zero. But then since $E\equiv 0/Z'$, 
we arrive at the situation in which $f(E)$ has codimension at least two. 
Now if $E\neq 0$, we get a contradiction with the negativity lemma: 
indeed by cutting by hypersurface sections, we can find a normal subvariety $Y'$ of $X'$ 
such that $Y'\to Z'$ is generically finite; now apply the negativity lemma 
to $E|_{Y'}$ and $-E|_{Y'}$ which are exceptional$/Z'$ to deduce that 
$E|_{Y'}=0$. So $E=0$ and $\phi^* L\sim_\R f'^*D$ for some $\R$-Cartier divisor $D$ on $Z'$.

If $Z$ is $\Q$-factorial, then by taking $R=\psi_*D$, we get $f^*R\sim_\R L$. We can then 
forget about $f'$ and rename $R$ to $D$.\\
\end{proof}

\section{Kodaira dimension of log divisors with big boundary}\label{s-Kodaira}

\subsection{Nef reduction maps of log divisors}
Theorem \ref{t-nef-dimension} demonstrates that nef reduction maps of log divisors are special. 
In dimension $3$ we can say much more (also see Proposition \ref{p-n=k} below). 

\begin{prop}\label{p-nef-red-log-div}
Let $(X,B)$ be a $\Q$-factorial projective klt pair of dimension $3$ over an uncountable 
$k$ of char $p>0$. Assume that $B$ is a big $\Q$-boundary and that $K_X+B$ is nef but 
not numerically trivial. 
Then there exist a  projective birational morphism $\phi\colon W\to X$ from a normal 
variety and a projective contraction $h\colon W\to Z$ to a normal variety, and a $\Q$-Cartier 
divisor $D$ on $Z$ such that 

$\bullet$ $\dim Z=n(K_X+B)$,

$\bullet$ $D$ is nef, and 

$\bullet$ $\phi^*(K_X+B)\sim_\Q h^*D$.
\end{prop}
\begin{proof}
We can assume  the nef dimension $n(K_X+B)=1$ or $2$ otherwise the statement is trivial.
Let $f\colon X\bir Z$ be a nef reduction map of $K_X+B$.

First assume $n(K_X+B)=2$. Then $Z$ is a surface and the singular locus of $X$ is vertical$/Z$, 
hence $X$ is smooth near the general fibres of $f$. Let $F$ be the generic fibre. 
Since $B$ is big, $B|_F$ is ample, and since $(K_X+B)|_F\equiv 0$, $K_F=K_X|_F$ is anti-ample. This implies that $\overline{F}$ is isomorphic 
to $\PP^1$ by [\ref{CTX}, Lemma 6.5] where $\overline{F}$ is the geometric generic fibre. 
Thus the pullback of $K_X+B$ to $\overline{F}$ is torsion and this is turn implies 
that $(K_X+B)|_F$ is torsion because $F$ is an integral scheme and $h^0(m(K_X+B)|_F)>0$ 
for some sufficiently divisible natural number $m$. Let $\phi\colon W\to X$ be a projective 
birational morphism from a normal variety such that the induced map $h\colon W\bir Z$ is a morphism. 
Now applying Lemma \ref{l-linear-pullback} to $h$ and $\phi^*(K_X+B)$, we can replace $W$ and $Z$ so that 
$\phi^*(K_X+B)\sim_\Q h^*D$ for some $\Q$-Cartier divisor $D$ on $Z$ which is necessarily nef.

We can then assume $n(K_X+B)=1$. Then $f$ is a morphism as it is regular over the generic 
point of $Z$ and $Z$ is a smooth curve. By Lemma \ref{l-relatively-trivial}, $K_X+B\equiv 0/Z$,  
and by Lemma \ref{l-num-pullback}, $K_X+B\equiv f^*D'$ for some $\Q$-divisor $D'$. Let $P=K_X+B-f^*D'$. 
We use an argument similar to [\ref{CTX}, proof of Theorem 1.9] to continue. 
Let $a\colon X\to \mathrm{Alb}$ be the Albanese morphism where $\mathrm{Alb}$ is the dual 
abelian variety of $\Pic^0(X)_{\rm red}$. Then for some sufficiently divisible 
natural number $m$, the divisor $mP$ belongs to $\Pic^0(X)_{\rm red}$ and $mP$ is the pullback of a numerically 
trivial divisor on $\mathrm{Alb}$. Now let $g\colon X\to S$ be the Stein factorization 
of $(f,a)\colon X\to Z\times \mathrm{Alb}$.  Assume $\dim S=1$. Then the induced morphism $e\colon S\to Z$ 
is an isomorphism which in turn implies that $a$ factors through $f$ and that 
$P$ is the pullback of some numerically trivial 
$\Q$-divisor $Q$ on $Z$, hence $K_X+B= f^*D$ where $D=D'+Q$. So we can take $W=X$ and $h=f$. 

Finally, we show that in fact $\dim S=1$. Assume otherwise. The arguments of the third paragraph of the proof of 
Lemma \ref{l-nef-thresh} show that every component of every fibre of $f$ over a closed point 
is uniruled. As $\mathrm{Alb}$ is an abelian variety, $a$ contracts each such component to a point or to a curve, hence $\dim S=2$. 
Let $A$ be an ample divisor on $S$.  As $B$ is big, we can assume that $B':=B-\epsilon g^*A\ge 0$ for some $\epsilon>0$. 
Since 
$$
K_X+B'\equiv -\epsilon g^*A/Z
$$
there is a $K_X+B'$-negative extremal ray $R$ over $Z$. 
By Theorem \ref{t-cone}, $R$ is generated by some rational curve $\Gamma$ 
which is contracted over $Z$. By construction, $\Gamma$ is not contracted over $S$ since $K_X+B'\equiv 0/S$. 
On the other hand, $\Gamma$ is contracted by $a$. 
Therefore, $\Gamma$ should be contracted over $S$ too, a contradiction.\\ 
\end{proof}

\subsection{ACC for horizontal coefficients}
The following result is of independent interest but also crucial for the proof of Proposition \ref{p-n=k} below.

\begin{prop}\label{p-ACC-global}
Let $\Lambda\subset [0,1]$ be a DCC set of real numbers. Then there is a 
finite subset $\Lambda^0\subset \Lambda$ with the following property: 
let $(X,B)$ be a pair and $f\colon X\to Z$ a projective contraction such that 

$\bullet$ $(X,B)$ is $\Q$-factorial dlt of dimension $3$ over $k$ of char $p>5$,

$\bullet$ $K_X+B$ is nef$/Z$ but not big$/Z$, 

$\bullet$ $B=\lambda H+B'$ where $H,B'\ge 0$, $H$ is big$/Z$, and $\lambda\in \Lambda$,

$\bullet$ the horizontal$/Z$ coefficients of $H$ and $B'$ are in $\Lambda$, and 

$\bullet$ $\dim X>\dim Z\ge 1$.\\\\
Then $\lambda$ belongs to $\Lambda^0$. 
\end{prop}
\begin{proof}
Assume $\dim Z=2$, let $F$ be the generic fibre of $f$, and $\overline{F}$ the geometric 
generic fibre. Since 
$K_X+B$ is not big$/Z$,  $(K_X+B)|_F\equiv 0$. As in the proof of Proposition \ref{p-nef-red-log-div}, 
one shows that $\overline{F}$ is isomorphic to $\PP^1$ and that $K_{\overline{F}}$ is the 
pullback of $K_X$. Let $B_{\overline{F}}$ be the pullback of $B$ to $\overline{F}$. 
Since $X$ is smooth near $F$, 
each coefficient of $B_{\overline{F}}$ is of the form 
$nb$ for some $n\in \N$ and some coefficient $b$ of $B$. In particular, the coefficients of $B_{\overline{F}}$
belong to some DCC set depending only on $\Lambda$. Therefore, these coefficients belong to some 
finite set depending only on $\Lambda$ because $\deg B_{\overline{F}}=2$. This in turn implies that 
the horizontal$/Z$ coefficients of $B$ belong to some finite set depending only on $\Lambda$, hence $\lambda$ belongs 
to some finite set depending only on $\Lambda$.

We can then assume $\dim Z=1$. Fix a fibre $F$ over a closed point $z$ 
such that $F$ and $B$ have no common components. 
Let $\phi\colon W\to X$ be a log resolution of $(X,B+F)$ so that over $U:=X\setminus \Supp F$, 
$\phi$ does not contract divisors with log discrepancy $0$ with respect to $(X,B)$. Such $\phi$ 
exist by the dlt assumption of $(X,B)$. 
Let $\Delta_W$ be the sum of the birational transform of $B$, plus the sum of the 
birational transform of the irreducible components of $F$ (each with coefficient $1$), 
and the sum of the exceptional prime divisors of $\phi$. Then $(W,\Delta_W)$ is dlt and  
 $\Supp \phi^*F\subset \rddown{\Delta_W}$. Moreover,   
if we write 
$$
K_W+\Delta_W=\phi^*(K_X+B)+E
$$
then we argue that $\Supp E\subseteq \rddown{\Delta_W}$: the exceptional$/X$ components of $E$ are 
obviously components of $\rddown{\Delta_W}$; on the other hand, 
$$
\phi_*\Delta_W-B=\phi_*E ~~\mbox{and}~~ \Supp(\phi_*\Delta_W-B)=\Supp F
$$ 
hence $\Supp \phi_*E=\Supp F$ which shows that the 
nonexceptional$/X$ components of $E$ are also components of $\rddown{\Delta_W}$. In addition, 
over $U$ the divisor $E$ is effective and its support is the reduced exceptional 
divisor of $\phi$ because $(X,B)$ is dlt, because over $U$ the boundary $\Delta_W$ 
is the sum of the birational transform of $B$ and 
the reduced exceptional divisor of $W\to X$, and because of our choice of $\phi$.

Run an LMMP$/X$ on $K_W+\Delta_W$ with scaling of some ample divisor as in [\ref{B}, 3.5]. 
By special termination [\ref{B}, Proposition 5.5], the LMMP terminates with a model $Y/X$ 
because it is an LMMP on $E$ 
and $\Supp E\subseteq \rddown{\Delta_W}$. 
In particular, since over $U$ the divisor $E$ is effective with support equal to the 
reduced exceptional divisor of $\phi$, the LMMP contracts any 
component of $E$ whose generic point maps into $U$. Thus $Y\to X$ is an 
isomorphism over $U$, and 
$E_Y$ maps into $\Supp F$. In particular, this means that 
$E_Y$ is supported on the fibre of the induced morphism $g\colon Y\to Z$ over the point $z=f(F)$.  

Now run an LMMP$/Z$ on $K_Y+\Delta_Y$ with scaling of some ample divisor. 
By special termination, the LMMP terminates near $\rddown{\Delta_Y}$. 
Since $K_Y+\Delta_Y$ is nef over $Z\setminus \{z\}$, the LMMP contracts only extremal rays over $z$ 
so such rays intersect the fibre over $z$, hence the LMMP terminates as the support of the 
fibre is inside $\rddown{\Delta_Y}$. Denote the resulting model 
by $Y'$, and denote the birational transform of $H$ by $H_{Y'}$. 
Then we can write $\Delta_{Y'}=\lambda H_{Y'}+\Delta_{Y'}'$, and  
by replacing $X$ with $Y'$, $H$ with $H_{Y'}$, and $B'$ with $\Delta_{Y'}'$, 
we can assume that $\Supp F\subseteq \rddown{B'}$ (we also need to add $1$ to $\Lambda$ 
since some of the coefficients of the new $B'$ are equal to $1$).

Let $S$ be a component of $F$ which intersects $H$. Note that $S$ is automatically normal as $(X,B)$ 
is $\Q$-factorial dlt [\ref{HX}][\ref{B}, Lemma 5.2]. By adjunction, we can write 
$$
K_S+B_S'=(K_X+B')|_S ~~\mbox{and}~~ K_S+B_S=(K_X+B)|_S 
$$
where $B_S:=\lambda H|_S+B_S'$ and the coefficients of $B_S$ and $B_S'$ belong to a DCC set determined by $\Lambda$ [\ref{B}, Proposition 4.2]. 
More precisely, the coefficient of each prime divisor $V\subset S$ in $B_S$ is of the form 
$$
\frac{l-1}{l}+\sum\frac{ b_i'm_i}{l}+\lambda \sum \frac{h_jn_j}{l}
$$ 
for some $l\in\N$ and $m_i,n_j\in\N\cup\{0\}$ where $b_i'$ and $h_j$ are the coefficients of 
$B'-S$ and $H$ respectively ($m_i>0$  if the corresponding component of $B'-S$ contains $V$; 
similarly $n_j>0$ if the corresponding component of $H$ contains $V$). 

We show that we can choose $S$ so that $H|_S$ is big. 
Since $H$ is big$/Z$, we can write $H\sim_\R A+N/Z$ where $A$ is ample and 
$N\ge 0$. Let $t$ be the smallest real number such that $N+tF\ge 0$. 
Then there is a component $S$ of $F$ which is not a component of $N+tF$.
But then $H|_S\sim_\R A|_S+(N+tF)|_S$ is big.  

Now since  $K_X+B$ is not big$/Z$, $K_X+B$ restricted to any fibre 
of $f$ is not big, so $(K_X+B)|_F$ is not big. This in turn implies that $K_S+B_S$ is not 
big. However, $K_S+B_S$ is semi-ample and it defines a contraction $S\to T$. 
Since $H|_S$ is big, $H|_S$ is horizontal$/T$. 
Applying [\ref{B}, Proposition 11.7], the horizontal$/T$ coefficients of $B_S$ belong to a 
finite set depending only on $\Lambda$. In particular, $\lambda$ belongs to a finite set 
depending only on $\Lambda$.\\ 
\end{proof}

\subsection{Kodaira dimension}

We come to the main result of this section.

\begin{prop}\label{p-n=k}
Let $(X,B)$ be a $\Q$-factorial projective klt pair of dimension $3$ over an uncountable
$k$ of char $p>5$. Assume that $B$ is a big $\Q$-boundary and that $K_X+B$ is nef but not numerically trivial. 
Then there exist a  projective birational morphism $\phi\colon W\to X$ from a normal 
variety and a projective contraction $h\colon W\to Z$ to a normal variety, and a $\Q$-Cartier 
divisor $D$ on $Z$ such that 

$\bullet$ $D$ is nef and big, and  

$\bullet$ $\phi^*(K_X+B)\sim_\Q h^*D$.\\\\ 
In particular, 
$$
\kappa(K_X+B)=\nu(K_X+B)=n(K_X+B)=\dim Z
$$
\end{prop}

\begin{proof} 
\emph{Step 1.} 
Assume there is a birational map $X\bir X'$ whose inverse does not contract any divisor 
and such that $(X',B')$ is $\Q$-factorial klt and the pullbacks of $K_X+B$ and $K_{X'}+B'$ 
are equal on some common resolution of $X$ and $X'$. If $X\bir X'$ contracts divisors, we replace $(X,B)$ with 
$(X',B')$. Repeating this process, we can assume that any map $X\bir X'$ as above is an isomorphism in 
codimension one.

Now if $K_X+B$ is big, then the proposition is trivial. 
So we can assume the numerical dimension $\nu(K_X+B)=1$ or $2$.
Moreover, by Theorem \ref{t-nef-dimension},  we can assume the nef dimension $n(K_X+B)$ is $1$ or $2$.  
By Proposition \ref{p-nef-red-log-div}, there exist a  projective birational morphism $\phi\colon W\to X$ from a normal 
variety and a projective contraction $h\colon W\to Z$ to a normal variety, and a $\Q$-Cartier 
divisor $D$ on $Z$ such that $\dim Z=n(K_X+B)$, $D$ is nef, and $\phi^*(K_X+B)\sim_\Q h^*D$. 
Moreover, the induced map $f\colon X\bir Z$ is a nef reduction map of $K_X+B$, 
so in particular, it is regular and proper over some nonempty open subset of $Z$.

If $n(K_X+B)=1$, then $D$ is ample and we are done. 
 Thus from now on we may 
assume $n(K_X+B)=2$, hence $\dim Z=2$. It remains to show that $D$ is big.
Suppose that $D$ is not big. We will derive a contradiction in several steps.\\  

\emph{Step 2.}  
Let $H_Z\ge 0$ be an ample $\Q$-divisor on $Z$, and let $H=\phi_*h^*H_Z$.
Since $B$ is big, perhaps after replacing $B,H$, we can assume $B\ge H+A$ where $A\ge 0$ is ample. 
Now let $\epsilon>0$ be a sufficiently small rational number. Then $K_X+B-\epsilon H$ is not pseudo-effective 
because by Lemma \ref{l-nef-pushdown}, 
$$
\phi^*(K_X+B-\epsilon H)\le \phi^*(K_X+B)-\epsilon h^*H_Z\sim_\Q h^*(D-\epsilon H_Z)
$$ 
and because $D-\epsilon H_Z$ is not pseudo-effective. Let $\delta$ 
be the smallest number so that $K_X+B-\epsilon H+\delta A$ is pseudo-effective. 
Note as $\epsilon$ is sufficiently small, $\delta$ is sufficiently small too. 
We can assume $(X,B-\epsilon H+\delta A)$ is klt. By Theorem \ref{t-mmodel},  
the pair has a log minimal model $(Y,B_Y-\epsilon H_Y+\delta A_Y)$, and by Theorem \ref{t-nef-dimension}, 
the nef dimension of $K_Y+B_Y-\epsilon H_Y+\delta A_Y$ is at most $2$.

We want to show that $X\bir Y$ is an isomorphism in codimension one. 
Run an LMMP on $K_X+B-\epsilon H+\delta A$ with scaling of some large multiple of 
$A$. Denote the steps of the LMMP by $X_i\bir X_{i+1}/Z_i$ which is either a flip 
or a divisorial contraction with $X_{i+1}=Z_i$. Assume that for each $i<l$, $X_i\bir X_{i+1}/Z_i$ 
is a flip and that it is $K_X+B$-trivial, i.e. $K_{X_i}+B_i$ is numerically trivial over $Z_i$. 
By Proposition \ref{p-polytope-rays} (5), $K_{X_l}+B_l$ is numerically trivial over $Z_l$,  
and by the first paragraph of Step 1, $X_l\to Z_l$ is a flipping contraction. 
Therefore, the LMMP is $K_X+B$-trivial and it does not contract any divisor. 
Moreover, the LMMP terminates with a log minimal model by Proposition \ref{p-term-2}. 
Since different log minimal models are isomorphic in codimension one, we deduce 
$X\bir Y$ is also an isomorphism in codimension one.\\

\emph{Step 3.}
Let $g\colon Y\bir V$ be a nef reduction map of $K_Y+B_Y-\epsilon H_Y+\delta A_Y$. 
By Step 2, we can assume $\dim V>0$ otherwise $K_Y+B_Y-\epsilon H_Y+\delta A_Y\equiv 0$, 
so $K_X+B-\epsilon H+\delta A\equiv 0$, and taking the limit when $\epsilon$ approaches $0$ we end up with 
$K_X+B\equiv 0$, a contradiction.

Let $\tau$ be the largest number such that $K_Y+B_Y-\tau H_Y$ is pseudo-effective 
over the generic point of $V$. Since $K_X+B$ is pseudo-effective, $K_Y+B_Y$ is pseudo-effective, 
so $\tau\ge 0$. On the other hand, since $A_Y$ is big, $K_Y+B_Y-\epsilon H_Y$ is not pseudo-effective over the generic point 
of $V$, hence $\tau<\epsilon$. 
We want to show that $\tau=0$. Assume not. 
 By construction, $(Y,B_Y-\epsilon H_Y)$ is klt. Moreover, by ACC for 
 lc thresholds [\ref{B}, Theorem 1.10], $(Y,B_Y)$ is lc because $\tau$ is sufficiently small 
 as $\epsilon$ is sufficiently small, and because we can write $B_Y-\tau H_Y=B_Y'+(1-\tau)H_Y$ where   
 the coefficients of $B_Y'$ and $H_Y$ belong to a fixed finite 
set depending only on $(X,B)$ and $H$. Thus  $(Y,B_Y-\tau H_Y)$ is klt 
because $\tau>0$ and $(Y,B_Y-\epsilon H_Y)$ is klt. 

Let 
$(Y',B_{Y'}-\tau H_{Y'})$ be a log minimal model of $(Y,B_Y-\tau H_Y)$ over some nonempty 
open subset $V'$ of $V$. By definition of $\tau$, $K_{Y'}+B_{Y'}-\tau H_{Y'}$ is not big$/V'$.   
On the other hand, $K_Y+B_Y+\delta A_Y\equiv \epsilon H_Y$ over the generic point of $V'$, hence $H_Y$ is 
big over $V'$ which in turn implies that $H_{Y'}$ is big over $V'$. 
So Proposition \ref{p-ACC-global} implies that $1-\tau=1$ or that $1-\tau$ is bounded away from $1$.
In our case $1-\tau=1$ is the only possibility because 
 $\tau$ is sufficiently small, hence $\tau=0$, a contradiction.\\

\emph{Step 4.}
Assume that $\dim V=2$.
By the last step, $K_Y+B_Y$ is pseudo-effective but not big over the generic point of $V$. 
So $K_Y+B_Y$ restricted to the geometric generic fibre of $g\colon Y\bir V$ is torsion as 
the geometric generic fibre is isomorphic to $\PP^1$ (cf. proof of Proposition \ref{p-nef-red-log-div}). 
Thus $K_Y+B_Y$ restricted to the generic fibre of $g$ is torsion too. This in turn implies 
$K_X+B\sim_\Q M$ for some $M$ whose support does not intersect the generic fibre of $g$.
Therefore, $(K_Y+B_Y)\cdot G=0$ for the general fibres $G$ of $g$.

Let $\alpha$ denote the map $X\bir Y$. By Lemma \ref{l-open-isom},  
there is an open subset $U$ of $X$ such that $\alpha|_U$ is an isomorphism and such that 
the complement of $U_Y:=\alpha(U)$ in $Y$ has dimension at most $1$.
Since $\dim V=2$, $Y\setminus U_Y$ is vertical$/V$, hence $G\subset U_Y$. 
Therefore, if $G^\sim$ on $X$ is the birational transform of $G$, then 
$(K_X+B)\cdot G^\sim=0$. So if $G$ is very general, then $G^\sim$ is inside some  
 very general fibre of $f \colon X\bir Z$. 
In particular, $H\cdot G^\sim=0$ because $H$ is vertical$/Z$.
But then  $H_Y\cdot G=0$, hence 
$$
0=(K_Y+B_Y-\epsilon H_Y+\delta A_Y)\cdot G=\delta A_Y\cdot G=\delta A\cdot G^\sim
$$
which contradicts the assumption that $A$ is ample.\\

\emph{Step 5.}
Now assume $\dim V=1$. Then $g\colon Y\bir V$ is a morphism and  $K_Y+B_Y-\epsilon H_Y+\delta A_Y\equiv 0/V$ 
by Lemma \ref{l-relatively-trivial}. 
By Step 3,  $(Y,B_Y)$ is lc and $K_Y+B_Y$ is pseudo-effective but not big$/V$. 
We want to argue that $(Y,B_Y)$ has a weak lc model 
over $V$. Indeed since $(Y,B_Y)$ is lc and $(Y,B_Y-\epsilon H_Y+\delta A_Y)$ is klt 
and $K_Y+B_Y-\epsilon H_Y+\delta A_Y\equiv 0/V$, we have 
$$
K_Y+B_Y\equiv 2(K_Y+B_Y-\tilde{\epsilon} H_Y+\tilde{\delta} A_Y)/V
$$ 
where $\tilde{\epsilon}=\frac{\epsilon}{2}, \tilde{\delta}=\frac{\delta}{2}$ and that 
$(Y,B_Y-\tilde{\epsilon} H_Y+\tilde{\delta} A_Y)$ is klt. Thus 
$(Y,B_Y)$ has a $\Q$-factorial weak lc model $(Y',B_{Y'})$ over $V$ such that $Y'\bir Y$ does not contract divisors.

Now by Theorem \ref{t-cone}, $K_{Y'}+B_{Y'}+P_{Y'}$ is globally nef where $P_{Y'}$ is the 
pullback of a sufficiently ample divisor on $V$. Let $Y'\bir S$ be a nef reduction map of 
$K_{Y'}+B_{Y'}+P_{Y'}$.   
If $C$ is a very general fibre of $Y'\bir S$, then $(K_{Y'}+B_{Y'}+P_{Y'})|_C\equiv 0$ 
and by our choice of $P_{Y'}$ we can actually assume that $(K_{Y'}+B_{Y'})|_C\equiv 0$. 
Since $Y'\bir X$ does not contract divisors and 
since $K_X+B$ is nef, by Lemma \ref{l-nef-pushdown}, $(K_X+B)|_{C^\sim}\equiv 0$  
where $C^\sim$ is the birational transform of $C$. 

Assume $\dim S=1$ and let $C_Z$ on $Z$ be the image of $C^\sim$. 
Then $D|_{C_Z}\equiv 0$ where $D$ is as in Step 1. This shows that  
the nef dimension $n(D)<2$, a contradiction.

Now assume $\dim S=2$. 
Denote $X\bir Y'$ by $\alpha'$. Since $\alpha'^{-1}$ does not contract divisors, by Lemma \ref{l-open-isom},  
there is an open subset $U'\subseteq X$ such that $\alpha'|_{U'}$ is an isomorphism and 
$Y'\setminus U_{Y'}'$ is of dimension at most $1$ where $U_{Y'}'=\alpha'(U')$. 
Since $\dim S=2$, $Y'\setminus U_{Y'}'$ is vertical$/S$, hence the very general fibres $C$ of 
$Y'\bir S$ are inside $U_{Y'}'$. Thus $(K_X+B)\cdot {C^\sim}=0$, from which
we deduce that $C^\sim$ is inside the very general fibres of $X\bir Z$. 
In particular, $H\cdot {C^\sim}= 0$, hence $H_{Y'}\cdot {C}= 0$. But then 
$$
0=(K_{Y'}+B_{Y'}-\epsilon H_{Y'}+\delta A_{Y'})\cdot C=\delta A_{Y'}\cdot C=\delta A\cdot C^\sim
$$
which contradicts the assumption that $A$ is ample.\\ 
\end{proof}

\section{Some semi-ampleness criteria}\label{s-semi-ample}

\begin{lem}\label{l-s-ample-criterion}
Let $f\colon X\to Z$ be a surjective morphism from a normal projective variety to a normal 
projective surface, over $k$ of char $p>0$. 
Assume that $L$ is a nef $\Q$-divisor on $X$ such that 

$\bullet$ $L\sim_\Q f^*D$ for some nef and big $\Q$-divisor $D$ on $Z$, and

$\bullet$ $L|_{f^{-1}\mathbb{E}(D)}$ is semi-ample.\\

Then $L$ is semi-ample.   
\end{lem}
\begin{proof}
It is enough to show that $D$ is semi-ample.
By taking hyperplane sections we can find a normal closed subvariety $Y\subset X$ such that the induced map $Y\to Z$ is generically finite. By replacing $X$ with $Y$ and $L$ with $L|_Y$, we can assume that $f$ is generically finite. 
Take the Stein factorization $g\ \circ\ h \colon X\to T\to Z$ of $f$ with $h$ birational and $g$ finite. 
As $g$ is finite, $D$ is semi-ample iff $g^*D$ is semi-ample and $g^{-1}\E(D)=\E(g^*D)$. 
So by replacing $Z$ with $T$ and $D$ with $g^*D$, we can assume $f$ is birational.  
Since $X$ and $Z$ are surfaces, $\mathbb{E}(L)=f^{-1}\mathbb{E}(D)\cup S$ where 
$S$ is disjoint from $f^{-1}\mathbb{E}(D)$ and $S$ is contracted by $f$.  
Thus $L|_{\mathbb{E}(L)}$ is semi-ample because $L|_{f^{-1}\mathbb{E}(D)}$ is semi-ample 
by assumption and $L|_S\sim_\Q 0$. Now by Keel [\ref{Keel}], $L$ is semi-ample and so as $f$ is birational, $D$ is semi-ample too. \\ 
\end{proof}

\begin{lem}\label{ewm}
Let $X$ be a normal projective variety over an uncountable $k$ of char $p>0$.  Suppose 
$L$ is a nef $\Q$-divisor on $X$ with equal Kodaira and nef dimensions $\kappa(L)=n(L)\leq 2$. 
Then $L$ is endowed with a map $X\to V$ to a proper algebraic space 
$V$ of dimension equal to $\kappa(L)$.
\end{lem}
\begin{proof}
 Since $\kappa(L)\ge 0$, we can assume $L\ge 0$.
There is a nef reduction map $f\colon X\bir Z$ to some normal projective variety $Z$, where $n(L)=\dim Z$.  
Replacing $X$ we may assume $f$ is a morphism.
Since $L\ge 0$ and $L|_{F}\equiv 0$, $L|_F=0$ where $F$ is the generic fibre of $f$. 
Moreover, by Lemma \ref{l-linear-pullback}, perhaps after replacing $X,Z,$ we can assume $L\sim_\Q f^* D$ for some 
nef $\Q$-divisor $D$ on $Z$, in particular, 
$$
\kappa(D)=\kappa(L)=n(L)=n(D)=\dim Z
$$
It is enough to show that $D$ is endowed with a map.

Suppose first that $\kappa(D)=2$.
Then $Z$ is a surface and $\E(D)$ is a finite union of curves with $D|_{\E(D)}\equiv 0$, so $D|_{\E(D)}$ is endowed with the 
constant map to a point, hence by [\ref{Keel}, Theorem 1.9],  $D$ is endowed with a map $Z\to V$.
Now suppose $\kappa(D)=1$.  Then $Z$ is a smooth curve and $D$ is a big divisor on it, which is therefore ample, hence 
we take $V=Z$. Finally, if $\kappa(D)=0$, then $D\equiv 0$ is endowed with the constant map $Z\to V$ to a point. 

Note that the proof shows that when $\kappa(L)=0 $ or $1$, then $L$ is actually semi-ample and 
$X\to V$ is the projective contraction associated to $L$.\\
\end{proof}

\begin{lem}\label{l-s-ample-criterion-2}
Let $X$ be a normal projective variety  over $k$ of char $p>0$, and $L$ a 
nef $\Q$-divisor on $X$ with Kodaira dimension $0\le \kappa(L)\le 2$. Assume that $L$ is endowed with a map 
$f\colon X\to V$ onto a proper algebraic space of dimension $\kappa(L)$.  
Moreover assume that $L|_F\sim_\Q 0$ for every fibre $F$ of $f$.
Then $L$ is semi-ample.   
\end{lem}

\begin{proof}

Taking the Stein factorization of $f$ we can assume $f$ is a contraction. 
Since $\kappa(L)\ge 0$, we can assume $L\ge 0$.
After applying Chow's lemma and replacing $X$  we can assume $f$ factors as $X\to Z\to V$ where 
$h\colon X\to Z$ is a projective contraction and $Z\to V$ is birational.  In the same way we may also assume $Z$ is smooth.
Since $L\equiv 0/Z$ and $L\ge 0$, 
Lemma \ref{l-linear-pullback} gives us a nef $\Q$-divisor $D$ on $Z$ such that $L\lin_\Q h^*D$.
Note that since $D$ is endowed with the birational map $Z\to V$, we deduce that $D$ is nef and big 
and that  $\kappa(D)=\kappa(L)=\dim Z=\dim V$.

Assume first that $\kappa(L)=2$. Then $\dim Z=2$, $\mathbb{E}(D)$ is contracted by 
$Z\to V$, and 
$L|_{h^{-1}\mathbb{E}(D)}$ is semi-ample as $L$ is torsion on 
the fibres of $X\to V$ and each connected component of $h^{-1}\mathbb{E}(D)$ is contained in such a fibre. 
Now apply Lemma \ref{l-s-ample-criterion}.
On the other hand suppose $\kappa(L) = 1$. Then $\dim Z=1$ which implies that $D$ is ample 
and that $L$ is semi-ample. Finally if $\kappa(L)=0$, then $Z$ is a point and $L$ is torsion.\\
\end{proof}

\section{Good log minimal models}\label{s-bpf}

\begin{rem}\label{r-log-smooth}
Let $(X,B)$ be a projective klt pair of dimension $3$ over $k$. Assume that 
$B=B'+A$ where $B',A\ge 0$ are $\Q$-divisors and $A$ is ample. Let $\phi\colon W\to X$ be a 
log resolution of $(X,B)$. We can write 
$$
K_W+B_W'=\phi^*(K_X+B')+E'
$$  
where $B_W'$ is a $\Q$-boundary, $\phi_*B_W'=B'$, $(W,B_W')$ is klt, $E'\ge 0$ is exceptional$/X$ and its support contains all the 
exceptional divisors of $\phi$. Since $\phi$ is obtained by a sequence of blowups with 
smooth centres, there is an exceptional$/X$ divisor $G$ which is ample$/X$. By the 
negativity lemma, $G\le 0$. Now $\phi^* A+\epsilon G$ is ample for some small $\epsilon>0$. 
Pick a general $A_W\sim_\Q \phi^* A+\epsilon G$. 
Then  we have 
$$
K_W+B_W'+A_W\sim_\Q \phi^*(K_X+B'+A)+E
$$  
where $E:=E'+\epsilon G\ge 0$ is exceptional$/X$ and its support contains all the 
exceptional divisors of $\phi$. 

It is easy to see from the definitions and using the negativity lemma 
that any log minimal model (resp. weak lc model) of $(W,B_W'+A_W)$ 
is also a log minimal model (resp. weak lc model) of $(X,B'+A)$.\\ 
\end{rem}

\begin{lem}\label{l-dlt-s-ample}
Let $(X,B=B'+A)$ be a projective $\Q$-factorial dlt pair of dimension $3$ over $k$ of char $p>5$ 
such that $B',A\ge 0$ are $\Q$-divisors, $A$ is ample, and $\rddown{B}=\rddown{B'}$. 
Assume $(Y,B_Y)$ is a weak lc model 
of $(X,B)$ such that $(Y,B_Y)$ is $\Q$-factorial dlt and $Y\bir X$ does not contract divisors. 
Then $(K_Y+B_Y)|_{\rddown{B_Y}}$ is semi-ample. 
\end{lem}
\begin{proof}
Since $(Y,B_Y)$ is dlt, its lc centres are the lc centres of $(Y,\rddown{B_Y})$, in particular, as
$\rddown{B_Y}=\rddown{B_{Y}'}$, $\Supp A_Y$ does not contain any lc centre.  
On the other hand, if $U\subseteq X$ is the largest open subset over which 
$\alpha\colon X\bir Y$ is an isomorphism, then $\Supp A_Y$ contains $Y\setminus \alpha(U)$ 
[\ref{B}, proof of Theorem 9.5]. So $Y\setminus \alpha(U)$ does not contain any lc centre 
of $(Y,B_Y)$.
 
Let $H_Y$ be a general ample $\Q$-divisor on $Y$ so that $A-H$ is also ample where $H$ is the 
birational transform of $H_Y$. Let $A'\sim_\Q A-H$ be general. Pick a small rational number $\epsilon>0$, 
and let 
$$
\Delta:=B'+(1-\epsilon) A+\epsilon A'+\epsilon H
$$ 
Then $(X,\Delta)$ is dlt on $U$, hence $(Y,\Delta_Y)$ is dlt on $\alpha(U)$. 
Moreover, since $(Y,B_Y)$ has no lc centre inside $Y\setminus \alpha(U)$, 
$\Supp (A_Y'+H_Y)$ does not contain any lc centre of $(Y,B_Y)$, so  
we see that $(Y,\Delta_Y)$ is dlt everywhere.

Now apply [\ref{B}, Theorem 1.9] to deduce that $(K_Y+\Delta_Y)|_{\rddown{\Delta_Y}}$ 
is semi-ample which in turn implies that $(K_Y+B_Y)|_{\rddown{B_Y}}$ is semi-ample 
because $K_Y+\Delta_Y\sim_\Q K_Y+B_Y$ and $\rddown{\Delta_Y}=\rddown{B_Y'}$.\\  
\end{proof}

\begin{prop}\label{p-fg-kappa=2}
Let $(X,B)$ be a projective  klt pair of dimension $3$ over $k$ of char $p>5$. 
Assume that $B$ is a big $\Q$-boundary and that $\kappa(K_X+B)\ge 0$. Then 
 $(X,B)$ has a good log minimal model. 
\end{prop}
\begin{proof}
We mimic the proof of [\ref{B}, Theorem 1.3]. We can assume $\kappa(K_X+B)\le 2$ by 
[\ref{B}, Theorems 1.2 and 1.4]. By extending $k$ we can assume it is uncountable.
Applying [\ref{B}, Lemma 7.7] we can assume $X$ is $\Q$-factorial.\\

\emph{Step 1.} 
By assumption $B$ is a big $\Q$-divisor, so by [\ref{Xu}][\ref{B}, Lemma 9.2], 
we can assume $B=B'+A$ where $B'$ is an effective $\Q$-divisor and $A\ge 0$ is an ample $\Q$-divisor.  
Take $M$ such that $K_X+B\sim_\Q M\ge 0$. By Remark \ref{r-log-smooth}, we can assume $(X,B+M)$ is log smooth.  
Replacing $A$ up to $\Q$-linear equivalence we may assume that $A$ and $B'+M$ share no components.
Since 
$$
K_X+B'+\epsilon M+A\sim_\Q (1+\epsilon)(K_X+B)
$$  
we may also replace $B'$ with $B'+\epsilon M$ 
for then we can assume $\Supp M\subseteq \Supp B'$.  

We want to prove  a slightly more 
general problem: we consider dlt pairs rather than klt. 
So we work with triples $(X,B,M)$ such that the following hold:
\begin{enumerate}
\item{$(X,B)$ is a dlt threefold pair over $k$ of char $p>5$, and $B$ is a $\Q$-boundary,}
\item{$K_X+B\sim_\Q M\ge 0$ for some $\Q$-divisor $M$,}
\item{$(X,B+M)$ is log smooth,}
\item{$0\leq \kappa(K_X+B)\leq 2$,}
\item{$B=B'+A$ where $B',A\ge 0$ are $\Q$-divisors, and $A$ is  ample with no common components with $B'+M$, and}
\item{$\Supp \rddown{B}\subseteq \Supp M \subseteq\Supp B'$.}\\
\end{enumerate}

Pick $(X,B,M)$ satisfying the above conditions. We will show that $(X,B)$ has a  
good log minimal model $(Y,B_Y)$ such that $Y\bir X$ does not contract divisors. 
Define $\theta(X,B,M)$ to be the number of components of $M$ which are not components of $\rddown{B}=\rddown{B'}$.\\

\emph{Step 2.} 
First assume $\theta(X,B,M)=0$, which is equivalent to $\Supp{M}\subseteq\rddown{B}$. 
We can run an LMMP on $K_X+B$ which terminates with a log minimal model $(Y,B_Y)$ by 
special termination [\ref{B}, Proposition 5.5]. By assumption, $\kappa(K_Y+B_Y)\ge 0$. 
If $K_Y+B_Y$ is numerically trivial, then it is torsion and we are done. If $K_Y+B_Y$ is not numerically trivial, 
then $\kappa(K_Y+B_Y)=n(K_Y+B_Y)$ by Proposition \ref{p-n=k}, hence 
$K_Y+B_Y$ is endowed with a map $Y\to V$ by Lemma \ref{ewm}. 
Moreover, $(K_Y+B_Y)|_{\rddown{B_Y}}$ is semi-ample by Lemma \ref{l-dlt-s-ample}.

Since $Y\to V$ is the map associated to $K_Y+B_Y$, on any of its fibres $F$ we have $M_Y|_F\num 0$. 
So either $F \cap \Supp M_Y=\emptyset$ or $F_{\rm red}\subseteq \Supp M_Y$.
 In the former case we have $M_Y|_F= 0$. In the latter case we show $M_Y|_F\sim_\Q 0$:  
from $F_{\rm red}\subseteq\Supp M_Y \subseteq \rddown{B_Y}$ and 
$(K_Y+B_Y)|_{\rddown{B_Y}}\lin_\Q M_Y|_{\rddown{B_Y}}$ being semi-ample 
we deduce $M_Y|_{F_{\rm red}}$ is semi-ample which in turn implies that 
$M_Y|_F$ is semi-ample  as char $k>0$; since $M_Y|_F\num 0$ we must have $M_Y|_F\lin_\Q 0$.  
Now by Lemma \ref{l-s-ample-criterion-2}, $M_Y$ is semi-ample. 
So from now on we can assume $\theta(X,B,M)>0$.\\

\emph{Step 3.}
Define 
$$
\alpha:=\min\{t>0|\rddown{(B+tM)^{\leq 1}}\neq \rddown{B}\}
$$
 where $D^{\leq1}$ is $D$ with all coefficients truncated at 1.
We can write $(B+\alpha M)^{\leq 1}=B+C$ and write $\alpha M = C+N$.  By construction $C,N\ge 0$ 
and  $\Supp N=\rddown{B}$ by (6).  By definition of $\alpha$ there must be some components of $\rddown{B+C}$ which are not in $\rddown{B}$, and these components are in $\Supp M$.  
The construction ensures that 
$$
\theta(X,B+C,M+C)<\theta(X,B,M)
$$

It is easy to show that $(X,B+C,M+C)$ satisfies the properties (1)-(6) of Step 1.
Indeed properties (1) and (3) follow from the assumption that $(X,B+M)$ is log smooth, 
and that $C$ is supported on $\Supp M$ and $B+C$ has coefficients at most 1. Properties 
(2),(5) and (6) are obvious. Property (4) is a consequence of 
$$
K_X+B\leq K_X+B+tC\leq K_X+B+\alpha M\lin_\Q (1+\alpha)(K_X+B)
$$ 
as it implies  
$$
\kappa(K_X+B) = \kappa(K_X+B+tC)
$$ 
for all rational numbers $t\in[0,1]$.

Let $\mathcal{T}$ be the set of those $t\in[0,1]$ such that $(X,B+tC)$ { has a good log minimal model $(Y,B_Y+tC_Y)$  
such that $Y\bir X$ does not contract divisors.} 
Arguing by induction on $\theta(X,B,M)$ and taking into account the previous paragraph, 
we can assume $1\in \mathcal{T}$ (note that the case $\theta=0$ of the 
induction was settled in Step 2).\\ 

\emph{Step 4.} 
Choose $0<t\in\mathcal{T}$.  We want to show that there is an $\epsilon>0$ such that $[t-\epsilon,t]\subset \mathcal{T}$.
Let $(Y,B_Y+tC_Y)$ be a good log minimal model of $(X,B+tC)$ such that $Y\bir X$ 
does not contract divisors. As $K_Y+B_Y+tC_Y$ is 
semi-ample, it defines a contraction $f\colon Y\to T$.  
Choose a sufficiently small $\epsilon>0$ and run a $K_Y+B_Y+(t-\epsilon)C_Y$-LMMP over $T$ with scaling of $\epsilon C_Y$ 
as in [\ref{B}, 3.5]. This is an LMMP on 
$$
M_{Y}+(t-\epsilon) C_{Y}=\frac{1}{\alpha} N_{Y}+(\frac{1}{\alpha}+t-\epsilon)C_{Y}
$$ 
and as $C_Y$ is positive on each extremal ray in the process, the LMMP is also an LMMP on $N_{Y}$. 
The LMMP terminates on some model $Y'$ by special termination [\ref{B}, Proposition 5.5] because $\Supp N_Y\subseteq \rddown{B_Y}$. 
Note that the LMMP is also an LMMP$/T$ on $K_Y+B_Y+(t-\epsilon')C_Y$ for any $\epsilon'\in (0,\epsilon)$. 
So we can replace $\epsilon$ with a smaller number if necessary. 
In particular, we can assume $t-\epsilon$ is rational and that $K_{Y'}+B_{Y'}+(t-\epsilon)C_{Y'}$ is globally nef 
by Theorem \ref{t-cone} or Proposition \ref{p-polytope-rays} (2).

If $K_{Y'}+B_{Y'}+(t-\epsilon) C_{Y'}$ is numerically trivial, then it is torsion as it has 
nonnegative Kodaira dimension. If it is not numerically trivial, then its Kodaira dimension and nef dimension are equal by 
Proposition \ref{p-n=k}, hence  $K_{Y'}+B_{Y'}+(t-\epsilon) C_{Y'}$ is 
endowed with a map $Y'\to V$ by Lemma \ref{ewm}. 
By construction, $K_{Y'}+B_{Y'}+t C_{Y'}$ is semi-ample and $\R$-linearly trivial over $T$. 
Moreover, $Y'\to T$ factors through $Y'\to V$ as $\epsilon$ is 
sufficiently small: indeed we can assume $K_{Y'}+B_{Y'}+(t-\delta) C_{Y'}$ 
is nef for some $\delta>\epsilon$, hence for any curve $\Gamma$ contracted by $Y'\to V$ 
we have $(K_{Y'}+B_{Y'}+(t-\epsilon) C_{Y'})\cdot \Gamma=0$ which implies that 
$(K_{Y'}+B_{Y'}+t C_{Y'})\cdot \Gamma=0$ and $(K_{Y'}+B_{Y'}+(t-\delta) C_{Y'})\cdot \Gamma=0$;
thus any 
curve contracted by $Y'\to V$ is also contracted by $Y'\to T$. 
We conclude that  
$$
K_{Y'}+B_{Y'}+tC_{Y'}\sim_\R 0/V
$$  
 In particular, as 
$K_{Y'}+B_{Y'}+(t-\epsilon)C_{Y'}\equiv 0/V$, we get $N_{Y'}\equiv 0/V$ and $C_{Y'}\equiv 0/V$.

Put $t':=t-\epsilon$. 
Let $F$ be a fibre of $Y'\to V$.  Since $N_{Y'}\equiv 0/V$, 
either $F \cap  \Supp N_{Y'}=\emptyset$ or $F_{\rm red} \subseteq  \Supp N_{Y'}$.  
In the first situation $N_{Y'}|_F=0$, so from $(K_{Y'}+B_{Y'}+tC_{Y'})|_F\sim_\R 0$ 
we deduce $C_{Y'}|_F\lin_\Q 0$ and $(K_{Y'}+B_{Y'}+t'C_{Y'})|_F\lin_\Q 0$.
In the second situation, by Lemma \ref{l-dlt-s-ample}, $(K_{Y'}+B_{Y'}+t'C_{Y'})|_{\rddown{B_{Y'}}}$
is semi-ample which implies that $(K_{Y'}+B_{Y'}+t'C_{Y'})|_F$ is semi-ample as 
$F_{\rm red} \subseteq  \Supp N_{Y'}\subseteq \rddown{B_{Y'}}$.
So 
$(K_{Y'}+B_{Y'}+t'C_{Y'})|_F\lin_\Q0$.  
Thus in any case the conditions of Lemma \ref{l-s-ample-criterion-2} are satisfied and so $K_{Y'}+B_{Y'}+t'C_{Y'}$ is semi-ample.
Therefore, $K_{Y'}+B_{Y'}+t''C_{Y'}$ is semi-ample and $({Y'},B_{Y'}+t''C_{Y'})$ is a good log minimal model of 
$(X,B+t''C)$ for any $t''\in [t-\epsilon,t]$, hence 
$[t-\epsilon,t]\subseteq \mathcal{T}$ as claimed.\\

\emph{Step 5.} 
Let $\tau:=\inf \mathcal{T}$. Assuming $\tau\notin \mathcal{T}$, we derive a contradiction. 
Take a strictly decreasing sequence of rational numbers $t_i\in\mathcal{T}$ approaching $\tau$.
For each $i$, there is a good log minimal model $(Y_i,B_{Y_i}+t_iC_{Y_i})$ 
of $(X,B+t_iC)$ such that $Y_i\bir X$ does not contract divisors. By taking a subsequence, we can assume that 
all the $Y_i$ are isomorphic in codimension one. In particular, 
$K_{Y_1}+B_{Y_1}+\tau C_{Y_1}$ is (numerically) a limit of movable divisors.
Run the LMMP on $K_{Y_1}+B_{Y_1}+\tau C_{Y_1}$ with scaling of $(t_1-\tau)C_{Y_1}$. 
Reasoning as in the first paragraph of Step 4, the LMMP terminates 
with a model $Y$ on which $K_{Y}+B_{Y}+\tau C_{Y}$ is nef. 
Note that the LMMP does not contract any divisor by the above movablity property. 
Moreover, $K_{Y}+B_{Y}+(\tau+\delta) C_{Y}$ is nef for some $\delta>0$. Now, by replacing the sequence, 
we can assume that $K_{Y}+B_{Y}+t_i C_{Y}$ is nef for every $i$ and by replacing each $Y_i$ with $Y$  
we can assume that $Y_i=Y$ for every $i$.  
A simple comparison of discrepancies (cf. [\ref{B-mmodel}, Claim 3.5]) shows 
that $(Y,B_{Y}+\tau C_{Y})$ is a $\Q$-factorial dlt 
weak lc model of $(X,B+\tau C)$. If we show that $K_Y+B_Y+\tau C_Y$ is semi-ample, then 
$\tau \in\mathcal{T}$ by Proposition \ref{p-term-2}, a contradiction.\\ 
 
\emph{Step 6.}
It remains to show that $K_Y+B_Y+\tau C_Y$ is semi-ample.
Let $Y\to T_i$ be the contraction defined by $K_{Y}+B_{Y}+t_i C_{Y}$. 
For each $i$, the map $T_{i+1}\bir T_i$ is a morphism because any curve 
 contracted by $Y\to T_{i+1}$ is also contracted by $Y\to T_i$. So perhaps after replacing the sequence, 
we can assume that $T_i$ is independent of $i$ so we can drop the subscript and simply 
use $T$. Note that $C_Y\sim_\Q 0/T$. 

Assume that $\tau$ is irrational. 
If $K_Y+B_Y+(\tau-\epsilon)C_Y$ is nef for 
some $\epsilon>0$, then $K_Y+B_Y+\tau C_Y$ is semi-ample because in this case 
$K_Y+B_Y+(\tau-\epsilon)C_Y$ is the pullback of a nef divisor on $T$ and 
$K_Y+B_Y+t_i C_Y$ is the pullback of an ample divisor on $T$. If there is no $\epsilon$ as above, then 
 there is a curve 
$\Gamma$ generating some extremal ray such that $(K_Y+B_Y+\tau C_Y)\cdot \Gamma=0$ 
and $C_Y\cdot \Gamma>0$ by [\ref{B}, 3.4] and Theorem \ref{t-cone}. 
This is not possible since $\tau$ is assumed to be irrarional. 
 So from now on we assume that 
$\tau$ is rational.
 
By Proposition \ref{p-n=k} and Lemma \ref{ewm}, $K_Y+B_Y+\tau C_Y$ is endowed with a map $f\colon Y\to V$. 
Any curve contracted by $Y\to T$ is also contracted by $Y\to V$. So
$Y\to V$ factors through $Y\to T$. 
Now $K_Y+B_Y+\tau C_Y\equiv 0/V$, so $C_Y$ is nef$/V$ but $N_Y$ is anti-nef$/V$. 
Let $F$ be a fibre of $Y\to V$. Then $F$ is disjoint from $\Supp N_Y$ or $F_{\rm red}$ is contained in $\Supp N_Y$. 
 Suppose $F \cap  \Supp N_Y=\emptyset$.  Then near $F$, 
$K_Y+B_Y+\tau C_Y$ is a positive multiple of $K_Y+B_Y+t_i C_Y$, hence 
$(K_Y+B_Y+\tau C_Y)|_F$ is semi-ample as $(K_Y+B_Y+t_i C_Y)|_F$ is semi-ample. 
Thus  $(K_Y+B_Y+\tau C_Y)|_F\sim_\Q 0 $.    
On the other hand, if $F_{\rm red}\subseteq  \Supp N_Y$, then again $(K_Y+B_Y+\tau C_Y)|_F\sim_\Q 0 $
because  $(K_Y+B_Y+\tau C_Y)|_{\rddown{B_Y}}$ is semi-ample by Lemma \ref{l-dlt-s-ample} 
and $F_{\rm red}\subseteq \Supp N_Y\subseteq \rddown{B_Y}$. 
The conditions of Lemma \ref{l-s-ample-criterion-2} are now satisfied, hence 
$K_Y+B_Y+\tau C_Y$ is semi-ample.\\
\end{proof}

\section{Proof of main results}\label{s-main-results}

Theorem \ref{t-cone} was already proved in Section \ref{s-cone}. We will give the proofs of 
the other main results.

\begin{proof}(of Theorem \ref{t-bpf})
We extend $k$ so that it is uncountable. By taking a $\Q$-factorialization [3, Lemma 6.7] 
we may assume that $X$ is $\Q$-factorial.  
Let $A=D-(K_X+B)$. By changing $A$ and $B$ we can assume $(X,\Delta:=B+A)$ is klt and that $A$ is an ample $\Q$-divisor. 
Moreover, if $P$ is the pullback of a sufficiently ample divisor on $Z$, then $K_X+\Delta+P$ 
is globally nef by Theorem \ref{t-cone}, and semi-ampleness of $K_X+\Delta+P$ implies semi-ampleness of 
$K_X+\Delta$ over $Z$. So replacing $\Delta$ with a boundary $\R$-linearly equivalent to $\Delta+P$, 
we can assume $Z$ is a point.

Assume that $D=K_X+\Delta\not\equiv 0$.   
By Proposition \ref{p-polytope-rays}, there exist real numbers $r_j$ and $\Q$-boundaries $\Delta_j$ 
such that $\Delta=\sum r_j \Delta_j$, $||\Delta-\Delta_j||$ 
are sufficiently small, $\Delta_j\ge A$, $({X},\Delta_j)$ are klt, and $K_X+\Delta_j$ are all nef. 
Moreover, we can ensure that $K_X+\Delta_j\not\equiv 0$ for every $j$.  By Proposition \ref{p-n=k}, 
$\kappa(K_X+\Delta_j)\ge 1$, hence $K_X+\Delta_j$ is semi-ample by Proposition \ref{p-fg-kappa=2}.
Therefore, $K_X+\Delta$ is semi-ample too. 

Now we can assume $D=K_X+\Delta\equiv 0$. As in the last paragraph, we can write 
$\Delta=\sum r_j \Delta_j$ making sure that $K_X+\Delta_j$ are numerically trivial 
$\Q$-divisors. This reduces the problem to the $\Q$-boundary case so we can assume 
$\Delta=B+A$ is already a $\Q$-boundary. Let $a\colon X\to \rm{Alb}$ be the albanese map 
and let $f\colon X\to T$ be the contraction given by the Stein factorization of $a$. 
Note that any rational curve on $X$ should  already be contracted by $f$ because 
$\rm{Alb}$ is an abelian variety. In particular, since $X$ is covered by rational curves 
(eg, see the end of the proof of Lemma \ref{l-nef-thresh}), we deduce that $f$ is not generically finite.

We will show that $T$ is a point. Assume not. 
Let $H$ on $X$ be the pullback of an ample $\Q$-divisor on $T$. Since $A$ is ample, 
can assume $A\ge H$. Pick a $K_X+B+A-H$-negative  
extremal ray $R$. By Theorem \ref{t-cone}, $R$ is generated by some rational curve $\Gamma$. 
Since $H\cdot R>0$, $\Gamma$ is not contracted by $f$, a contradiction. 
Therefore the albanese map is constant. But then some positive multiple $m(K_X+\Delta)$ is 
the pullback of a divisor on $\rm{Alb}$, hence $m(K_X+\Delta)\sim_\Q 0$.\\
\end{proof}

\begin{proof}(of Theorem \ref{t-contraction})
Let $R$ be a $K_X+B$-negative extremal ray$/Z$. Note that $R$ being over $Z$ means that 
$P\cdot R=0$ where $P$ is the pullback of some ample divisor on $Z$. By adding a small ample 
divisor to $B$ and by perturbing the coefficients, we can assume $B$ is a big $\Q$-boundary. 
In particular, there are only finitely many negative extremal rays of $K_X+B$ and  
 they are all generated by extremal curves with bounded intersection with $K_X+B$, by Theorem \ref{t-cone}. 
Thus there is an ample $\Q$-divisor $H$ such that $L=K_X+B+H$ is globally nef and $L^\perp=R$. 
By Theorem \ref{t-bpf}, $L$ is semi-ample so it defines a projective contraction $X\to T$.
The morphism $X\to T$ is nothing but the contraction of $R$. Since $R/Z$, the morphism $X\to Z$ 
factors through $X\to T$.\\
\end{proof}

\begin{proof}(of Theorem \ref{t-finiteness})
This follows from Theorem \ref{t-bpf} and Proposition \ref{p-finiteness}. 
\end{proof}

\begin{proof}(of Theorem \ref{t-term})
This follows from Theorems \ref{t-bpf} and \ref{t-contraction} and Proposition \ref{p-term}. 
\end{proof}

\begin{proof}(of Theorem \ref{t-term-2})
If $K_X+B$ is pseudo-effective$/Z$, this is already proved in Proposition \ref{p-term-2}. 
If $K_X+B$ is not pseudo-effective$/Z$, then the LMMP is also an LMMP on $K_X+B+\frac{1}{2}C$ 
with scaling of $\frac{1}{2}C$, hence it terminates by Theorem \ref{t-term}. 
\end{proof}

\begin{proof}(of Theorem \ref{t-Mfs})
We can find a projective contraction $\overline{f}\colon \overline{X}\to \overline{Z}$ of normal projective varieties 
such that $X$ is an open subset of $\overline{X}$ and $\overline{f}$ restricted to $X$ is $f$. 
Let $\phi\colon \overline{W}\to \overline{X}$ be a log resolution such that any prime exceptional 
divisor of $\phi$ whose generic point maps into $X$, has positive log discrepancy with respect 
to $(X,B)$. Let $B_{\overline{W}}$ be the 
sum of the birational transform of $B$ and the reduced exceptional divisor of $\phi$. Run an 
LMMP$/\overline{X}$ on $K_{\overline{W}}+B_{\overline{W}}$ with scaling of some ample divisor. 
By our choice of $\phi$ we reach a 
model $\overline{Y}$ such that $\overline{Y}\to \overline{X}$ is a small morphism 
over $X$. So we can replace $(X,B)$ with $(\overline{Y},B_{\overline{Y}})$, hence assume 
$X$ is projective and $\Q$-factorial. 
  
Pick a general sufficiently ample divisor $A$ so that $K_X+B+A$ is nef$/Z$ and $(X,B+A)$ is dlt. 
Let $\epsilon >0$ be small enough so that $K_X+B+\epsilon A$ is not pseudo-effective$/Z$. 
We can find a boundary $\Delta\sim_\R B+\epsilon A/Z$ such that $(X,\Delta+(1-\epsilon)A)$ is klt. 
Now run an LMMP$/Z$ on $K_X+\Delta$ with scaling of $(1-\epsilon)A$. By Theorem \ref{t-term}, 
the LMMP terminates with a Mori fibre space of $(X,\Delta)$ over $Z$. The LMMP is also an LMMP 
on $K_X+B$ with scaling of $A$, hence the Mori fibre space is also a Mori fibre space of $(X,B)$ over $Z$.\\
\end{proof}


\vspace{2cm}

\flushleft{DPMMS}, Centre for Mathematical Sciences,\\
Cambridge University,\\
Wilberforce Road,\\
Cambridge, CB3 0WB,\\
UK\\
c.birkar@dpmms.cam.ac.uk\\
jaw66@cam.ac.uk

\end{document}